\documentclass[11pt]{article}
\usepackage{amsfonts}
\usepackage{amsthm}
\usepackage{multirow}

\providecommand{\tabularnewline}{\\}

\usepackage{amsmath}
\usepackage{amssymb}
\usepackage{graphicx}
\usepackage{color}

\usepackage{marginnote}
\usepackage{amstext,amsthm,amssymb,amsmath}
\usepackage{graphicx}
\usepackage{subfigure}

\hfuzz=\maxdimen
\graphicspath{ {./pic/}}
\title{Multicontinuum homogenization and its relation to nonlocal multicontinuum theories}
\author{Yalchin Efendiev\footnote{Department of Mathematics, Texas A\&M University, College Station, TX 77843, USA}, ~ Wing Tat Leung\footnote{Department of Mathematics, City University of Hong Kong, Hong Kong}}

\begin{document}
\maketitle

\begin{abstract}
In this paper, we present a general derivation of multicontinuum
equations and discuss cell problems. We present 
constraint cell problem formulations in a representative 
volume element and oversampling techniques that allow reducing
boundary effects. We discuss different choices of constraints for cell
problems. We present numerical results that show how oversampling 
reduces boundary effects. Finally, we discuss the relation of the proposed
methods to our previously developed methods, 
Nonlocal Multicontinuum Approaches.

\end{abstract}

\section{Introduction}

One of commonly used approaches for multiscale problems includes
homogenization and its variations \cite{blp78,jko94,lipton06,ms00,weh02,bp04,dur91}, 
where effective properties at each macroscale
grid or point are computed. These computations are often based on
local solutions computed in a representative volume element
(or coarse grid)
centered at macroscale point. Homogenization-based approaches
assume scale separation and
 that the local media can be replaced by a homogeneous material.
As a result, it is assumed that the solution average in each 
coarse block approximates the heterogeneous solution within this 
coarse block.

In many cases, even within scale separation realm, 
homogenization (as discussed above) is not sufficient
and the coarse-grid formulation
 requires multiple homogenized coefficients.
These approaches are developed for different applications
\cite{rubin48,barenblatt1960basic,showalter1991micro,showalter1991micro,aifantis1979continuum,iecsan1997theory,bunoiu2019upscaling,arbogast1990derivation,bedford1972multi,chai2018efficient,alotaibi2022generalized} and
we call them 
 (following the literature)
multicontinuum approaches.
Multicontinuum approaches assume that the solution average is not
sufficient to represent the heterogeneous solution within
each coarse block.
In the derivation of multicontinuum approaches, there are typically
several assumptions: (1)  continua definitions; (2) physical laws 
describing the interaction among continua; and 
(3) conservation laws deriving final
equations. Various assumptions are typically made in deriving these models.
The first such approach is presented in
\cite{rubin48}, where the author assumes existence of continua that
have different equilibrium temperatures among each other (continua) and
formulates empirical laws for interaction among continua. 
In more rigorous approaches related to porous media
\cite{barenblatt1960basic,showalter1991micro,arbogast1990derivation}, 
the continua are assumed to be fracture and matrix regions.
In our earlier works \cite{chung2018nonlinear,chung2018non,zhao2020analysis}, 
we define 
the continua 
via local spectral decompositions and show that the resulting approach
converges independent of scales and contrast if representative volumes
are chosen to be coarse blocks. 
In this work, we use similar ideas (as in 
\cite{chung2018nonlinear,chung2018non,zhao2020analysis})
for problems with scale separation
and formulate cell problems and formally 
derive multicontinuum equations.

The main objectives of this paper are the following.
\begin{itemize}

\item We derive multicontinuum methods using a
homogenization-like expansion and
present constraint 
cell problem formulations.

\item Constraint cell problems allow using averages for different
quantities and regions (continua) and give flexibility to the framework.

\item We discuss appropriate local boundary conditions in 
representative volume elements for problems with scale separation
and introduce oversampling. Using oversampling, we consider
reduced constraint cell problems, where we use constraints
for the averages only. 

\item The resulting multicontinuum equations show that
local averages of the solution will differ among each other
if diffusion and reaction terms in the upscaled equations
 balance each other. This requires
smaller reaction and/or larger diffusion terms, which occur
in the presence of high contrast. We discuss this issue and 
show that a  multicontinuum concept is via local
spectral decomposition.

\item We discuss the relation to NLMC approaches that go beyond scale 
separation.

\item The average constraints, discussed in this paper,
 are easy to set and guarantee exponential decay (i.e., we remove boundary effects).

\item We present numerical results.

\end{itemize}

We note that 
to go beyond scale separation, numerical approaches use entire coarse blocks
(see Figure \ref{fig:ill}) to do local computations
\cite{chung2016adaptiveJCP,GMsFEM13,eh09,hw97,jennylt03,chung2018constraint,oz06_1}. Among these approaches, multiscale finite element method
and its variations are proposed, where multiscale basis functions
are computed on coarse grids.

\subsection{The main idea of this paper}

In this section, we briefly present the main idea of the 
paper. We consider a steady-state or dynamic problem
\[
\mathcal{L}(u) = f
\]
subject to some boundary and initial conditions. It is assumed
that the problem is solved on a computational grid consisting
of grid blocks, denoted $\omega$, that are much larger than
heterogeneities.
We assume some type of homogeneity within each
computational block represented by Representative 
Volume Element $R_\omega$ that corresponds to 
a computational element $\omega$ (see Figure \ref{fig:ill})
(more precise meaning will be defined later).
We assume that within each $R_\omega$, 
there are several distinct average states
can be achieved (known as multicontinua).
We denote the characteristics function for 
continuum $i$ within $R_\omega$ by $\psi_i^\omega$
($\omega$ will be omitted since local computations are restricted to a
coarse block), 
i.e., $\psi_i=1$ within
continuum $i$ (can be irregular shaped regions consisting of several
parts, in general) and $0$ otherwise.
We introduce oversampled $R_\omega^+$ that contains several $R_\omega^p$'s.
We denote the central (target) RVE by simply $R_\omega$.
In general, one can define the regions corresponding to the continuum
via local spectral decomposition of the solution space within $R_\omega$,
as discussed later.

We assume a variational formulation of the problem
\[
\sum_\omega \int_\omega \mathcal{A}(u,v)  = \int_\Omega fv,
\]
where $\mathcal{A}$ is the corresponding bilinear form.
We assume that $R_\omega$ can be used to approximate each integral
$\int_\omega$ (in general space-time integral). 
I.e.,
\begin{equation}
\label{eq:idea2}
 \int_\omega \mathcal{A}(u,v) \approx {|\omega|\over |R_\omega|}\int_{R_\omega} \mathcal{A}(u,v).
\end{equation}
Summation over repeated indices is assumed in the paper. 
Next, we construct local cell problems in $R_\omega$ that are used to represent
$u$.

We assume there are several macroscopic quantities 
denoted by $U_i^\omega$ in each $R_\omega$, 
where $i$ is the continuum. These quantities can be thought of
as average solutions within each continuum.
We introduce cell problems in $R_\omega^+$ 
(that consists of $R_\omega^p$)
that can distinguish 
these states. 
The first represents averages (formally written)
\begin{equation}
\begin{split}
\mathcal{L}(\phi_i) = r_i\ \ \text{in}\ R_\omega^+\\
\sum_p \int_{R_\omega^p} \phi_i \psi_j^p = \delta_{ij}\int_{R_\omega^p}  \psi_j^p ,
\end{split}
\end{equation}
where $r_i$ accounts for constraints,
and the second one accounts for the gradients (formally written)
\begin{equation}
\begin{split}
\mathcal{L}(\phi_i^m) = r_i^m\ \ \text{in}\ R_\omega^+\\
\int_{R_\omega^p} \phi_i^m \psi_j^p = \delta_{ij}\int_{R_\omega^p} (x_m-c) \psi_j, 
\ \ \\
\int_{R_\omega} (x_m-c) \psi_j^{p_0}=0.
\end{split}
\end{equation}
additional initial conditions are posed. $p_0$ refers to the target
RVE, $R_\omega$.
These cell problems are written formally and will more precisely be
described in next sections. We will use oversampling regions
and constraints in each $RVE$ within the oversampled region to avoid
boundary effects.
Using these cell problems, the local solution
in $R_\omega$ is written as
\begin{equation}
\label{eq:idea1}
u\approx \phi_i U_i + \phi_i^m \nabla_m U_i.
\end{equation}
We assume $U_i(x)$ is smooth function representing the $i$th continuum.
I.e., $U_i(x_\omega)\approx \int_{R_\omega} u \psi_i / \int_{R_\omega} \psi_i$,
with $x_\omega$ being a center point of $R_\omega$.
Substituting (\ref{eq:idea1}) into (\ref{eq:idea2})
and taking $v\approx \phi_s V_s + \phi_s^k \nabla_k V_s$,
we obtain multicontinuum equations for $U_i$.
Substituting (\ref{eq:idea1}) and the form for $v$ into 
equations, we obtain multicontinuum models.

Our main contributions are the following.

\begin{itemize}

\item We formulate constrained cell problems using Lagrange multipliers.

\item To avoid boundary effects, we formulate cell problems in oversampled regions and use solutions' averages to 
get fast decay of boundary effects. This is also shown numerically.

\item We derive multicontinuum upscaled models and
formulate scaling for each term, which is related to RVE size.
This
shows that unless there is some
type of high contrast, the averages $U_i$ within $R_\omega$ will become similar.

\item We formulate spectral continuum definitions, which
can be used to define $\psi_i$'s.

\item We discuss cell problems that use multiple constraints (averages and gradients) and discuss the advantages/disadvantages associated with such constraints.

\end{itemize}

The paper is organized as follows.
In the next section, we present preliminaries and show the arguments
used in \cite{rubin48}. Section 3 is devoted to the derivation of multicontinuum approaches for a scale separation case. In Section 4, we present
spectral continuum ideas. Section 5 is devoted to some remarks
that include the derivation using multiple constraints
and nonlinear multicontinuum models. Finally, we present some numerical
results in Section 6.

\section{Preliminaries}

\subsection{The work of L. I. Rubinstein \cite{rubin48} from 1948}

First, we briefly discuss the paper by L. I. Rubinstein \cite{rubin48},
which is the first paper that derives multicontinuum equations
based on physical laws. We skip/simplify some details. 
In \cite{rubin48}, the author considers time-dependent
diffusion equation in 
heterogeneous media. The equation
at the fine scale is
\begin{equation}
\begin{split}
u_t - \nabla\cdot(\kappa \nabla u) =f.
\end{split}
\end{equation}
The paper \cite{rubin48} assumes that the media consists of many small particles (possible
connected) divided into the group of  $N$ components (continua), where
the diffusivity of each component is $\kappa_i$.
The media is assumed to be stochastic, i.e., $\kappa(x,\zeta)$,
where $\zeta$ refers to a realization. 
At each point $x$,  $\omega_x$ is an elementary volume
around point $x$. We
denote the distribution within a component $i$
 as $\widetilde{U_i}(x,t,\zeta)$ and denote by
\[
U^*_i(x,t,\zeta)={1 \over \omega_x}\int_{\omega_x} \widetilde{U_i}(z,t,\zeta) dz
\]
and denote (mathematical expectation)
\[
U_i(x,t)=\int U^*_i(x,t,\zeta) d\nu(\zeta).
\]

It is assumed that within a representative element, different components 
can have 
different averages and  conservation for 
each component is written down.
 The conservation 
consists of three terms. The first term is the diffusion flux and is 
taken by (in \cite{rubin48})
\[
q_{1i}=\int_{\Sigma'} \kappa_i {\partial U_i\over\partial n} \mu_i d\sigma \Delta t,
\]
where $\mu_i$ is fraction of $i$th component on (larger) elementary volume 
boundary $\Sigma'$, $\Delta t$ is a time interval.
There are a number of assumption about components' homogeneities
 on boundaries of 
$\Omega'$ (RVE).
The second flux is taken to be heat exchange  within an elementary
volume, which occurs because of different temperatures within
each component. 
Using Henry's law, this flux 
is written in \cite{rubin48} as
\[
q_{2ij}=\int_{\Omega'} \alpha_{ij}^* (U_j - U_{ij}) d\omega \Delta t,
\]
where $U_{ij}$ is a temperature in $j$th component when $i$th component
temperature is $U_i$. It is taken to be $U_{ij}=U_i$.
The third flux is given by
\[
q_{3i} = \int_{\Omega'} c_i \rho_i {\partial U_i \over\partial t} p_i d\omega \Delta t,
\]
where $c_i$, $\rho_i$ represent fluid properties and $p_i$
is a volume fraction of $i$th component.
From
\[
q_{3i}=q_{1i} + q_{2i}
\]
one arrives to
\[
\nabla \cdot(\kappa_i\mu_i \nabla U_i) +  \sum_j \alpha_{ij}^*(U_j - U_{i}) = 
c_i \rho_i p_i {\partial U_i\over \partial t}.
\]
Similar multicontinuum models are proposed in different application
areas.

In this paper, we give a derivation based on formal
expansion, cell problems, and then show a relation to
theories developed in \cite{chung2018non,chung2018constraint}.
This derivation can be made rigorous under some assumptions
(cf. \cite{chung2021nonlocal}). We mention some assumptions as we go 
along without making them formal to keep the presentation simple.
Our derivation (1) reveals the nature of continua,
(2) shows their
relation to local spectral decomposition, and 
(3) formulate constraint cell problems with appropriate boundary
conditions.

\section{Multicontinuum derivation based on volume average constraints}

\subsection{Steady-state case}

\begin{figure}
\centering
\includegraphics[scale=0.35]{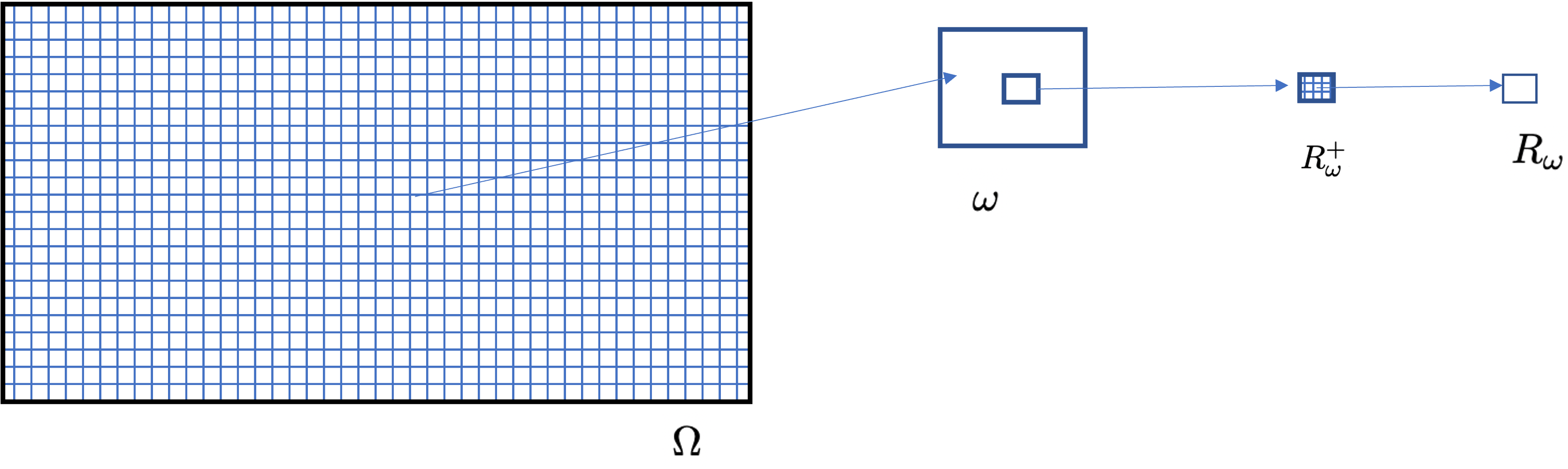}
\caption{Illustration}
\label{fig:ill}
\end{figure}

In this section will repeat some parts of Introduction.
Our approach starts from a finite element method formulated
on a coarse grid. Coarse grid contains RVE, where local
computations will be performed (see Figure \ref{fig:ill}).
We assume a partition of the domain into elements, where
 $\omega$ is a generic
coarse-grid element 
(triangle or rectangle), $R_\omega$ is a representative volume (RVE) within 
$\omega$ (see Figure \ref{fig:ill}). 
We consider a steady state diffusion equation
\[
\int_\Omega \kappa \nabla u \cdot \nabla v= \int_\Omega f v, \ v\in H_0^1(\Omega).
\]
Representative volume, as usual, is assumed to represent
the whole coarse block $\omega$ in terms of heterogeneities.
In each $R_\omega$, we assume $N$ continua (components) and introduce
\[
\psi_j=\text{ {1 in continuum $j$, 0 otherwise.}}
\]
In general, one can use different functions, e.g., eigenfunctions
of local problems \cite{chung2018non,chung2018constraint} to represent each continua, as discussed later.

Next, we remind that  $R_\omega^+$ is an oversampled 
region (RVE) that surrounds $R_\omega$. It is taken
to be several times larger compared to $R_\omega$
and consists of several RVE's, denoted by $R_\omega^p$ ($p$ is the numbering).
In general, they ($R_\omega^p$'s) 
can be the copies of $R_\omega$ and it is used
to remove boundary effects.
The target RVE, we denote by $R_\omega^{p_0}$ or simply $R_\omega$.
We
introduce 
two sets of cell problems with constraints.
\begin{equation}
\label{eq:cell_grad_aver}
\begin{split}
\int_{R_\omega^+}\kappa \nabla \phi^m_i\cdot \nabla v -  \sum_{j,p} {\beta_{ij}^{mp}\over \int_{R_\omega^p}\psi_j^p} \int_{R_\omega^p}\psi_j^p  v  =0\\
\int_{R_\omega^p}  \phi^m_i \psi_j^p = \delta_{ij} \int_{R_\omega^p} (x_m-c_{mj})\psi_j^p ,\\
\int_{R_\omega^{p_0}} (x_m-c_{mj})\psi_j^{p_0} =0\ \text{condition for $c$}, \\
\end{split}
\end{equation}
and
\begin{equation}
\label{eq:cell_aver_aver}
\begin{split}
\int_{R_\omega^+}\kappa \nabla \phi_i\cdot \nabla v -
\sum_{j,p} {\beta_{ij}^p\over \int_{R_\omega^p}\psi_j^p} \int_{R_\omega^p}\psi_j^p  v =0\\
\int_{R_\omega^p}  \phi_i \psi_j^p = \delta_{ij} \int_{R_\omega^p} \psi_j^p.
\end{split}
\end{equation}
The first cell problem (\ref{eq:cell_grad_aver})
one accounts for the gradient effects and is taken to vanish
in the target RVE, $R_\omega^{p_0}$. This cell problem accounts for
standard homogenization effects.
The second cell problem
(\ref{eq:cell_aver_aver}) accounts for different averages in
each continuum. By imposing the same averages in each $R_\omega^p$,
we reduce the boundary effects in an exponential manner
\cite{chung2018non,chung2018constraint}.
Here, for simplicity,
we do not use $\omega$ index in $\phi_i$ or $\phi_i^m$, though
both of them depend on $\omega$. This is because our calculations will
be done in each $\omega$ separately.
In general, one can remove the index $p$ in $\psi_j$ if RVE's are 
periodically repeated or similar.

Next, we formulate some properties of $\beta$'s. We note 
\begin{equation}
\begin{split}
\sum_{j,p} \beta_{ij}^p=0
\end{split}
\end{equation}
which can be obtained by taking $v=1$ in (\ref{eq:cell_aver}).
If we take $v=\phi_k$ in (\ref{eq:cell_grad_aver}), then,
we have
$\int_{R_\omega^+} \kappa \nabla \phi_i^m\cdot \nabla \phi_k = \sum_p \beta_{ik}^{mp}$.
If
 we take $v=\phi_s$ in (\ref{eq:cell_aver_aver}),
we get 
\begin{equation}
\label{eq:beta_aver}
\begin{split}
\beta_{is}^*=\int_{R_\omega^+} \kappa \nabla \phi_i\cdot \nabla \phi_s = \sum_p \beta_{is}^p.
\end{split}
\end{equation}

We assume that in $R_\omega$,
\begin{equation}
\label{eq:approx1}
\begin{split}
u\approx  \phi_i U_i  + \phi_i^m \nabla_m U_i,
\end{split}
\end{equation}
where $U_i$
is a smooth function representing the homogenized solution
for $i$th continuum. 
More precisely, $U_i$ can be thought as a limit 
of $\int_{R_\omega} u \psi_i / \int_{R_\omega}  \psi_i$ (piecewise constant function)
taken over all
$R_\omega$ as the RVE size goes to zero. 
We will assume $U_i$ and their gradients can be approximated by
 constants in RVE and use mid point to represent their values.
We note that (\ref{eq:approx1}) can be shown under the assumption that
$U_i$ is a smooth function. 

Next, we derive multicontinuum equations for $U_i$.
For any $v\in H_0^1$, we have
\begin{equation}
\label{eq:assump1}
\begin{split}
\int_\Omega f v = \int_\Omega \kappa \nabla u \cdot \nabla v = \sum_\omega \int_\omega \kappa \nabla u \cdot \nabla v \approx
\sum_\omega {|\omega|\over |R_\omega|} \int_{R_\omega} \kappa \nabla u \cdot \nabla v, 
\end{split}
\end{equation}
where we make an assumption that integrated average over
RVE can represent the whole 
computational element $\omega$. This approximation 
holds if all $ \int_{R_\omega^p}$ are approximately equal for all 
$R_\omega^p$ in $\omega$.
Next, we approximate each term 
\begin{equation}
\begin{split}
\int_{R_\omega} \kappa \nabla u \cdot \nabla v =
\int_{R_\omega} \kappa \nabla (\phi_i U_i) \cdot\nabla v  + 
\int_{R_\omega} \kappa \nabla (\phi_i^m \nabla_m U_i)   \cdot\nabla v. 
\end{split}
\end{equation}
We assume that the variation of $U_i$ and $\nabla_m U_i$
are small compared to
the variations of $\phi_i$ and $\phi_i^m$ (see scalings (\ref{eq:scalings})) 
and assume
$\int_{R_\omega} \kappa \nabla (\phi_i U_i) \cdot \nabla v\approx
\int_{R_\omega} \kappa (\nabla \phi_i) U_i \cdot\nabla v$
and
$\int_{R_\omega} \kappa \nabla (\phi_i^m \nabla_m U_i)  \cdot\nabla v
\approx
\int_{R_\omega} \kappa \nabla (\phi_i^m) \nabla_m U_i \cdot \nabla v
$.
We take $$v=\phi_s V_s + \phi_s^k\nabla_k V_s.$$
Then, denoting for simplicity, $R_\omega=R_\omega^{p_0}$, we
have
\begin{equation}
\begin{split}
\int_{R_\omega} \kappa \nabla (\phi_i U_i) \cdot \nabla v \approx
U_i(x_{\omega}) \int_{R_\omega} \kappa \nabla \phi_i \cdot \nabla v =\\
U_i(x_{\omega}) V_s(x_\omega)\int_{R_\omega} \kappa \nabla \phi_i \cdot \nabla \phi_s + U_i(x_{\omega}) \nabla_m V_s(x_\omega)\int_{R_\omega} \kappa \nabla \phi_i \cdot \nabla \phi_s^m =\\
 U_i(x_{\omega})  \beta_{is}^*   V_s(x_\omega) + \beta_{is}^{m*}  U_i(x_{\omega}) \nabla  V_s(x_\omega),
\end{split}
\end{equation}
where
\[
\beta_{ik}^{m*}=\int_{R_\omega} \kappa \nabla \phi_i^m \cdot \nabla \phi_k, \ \ \beta_{ik}^{*}=\int_{R_\omega} \kappa \nabla \phi_i \cdot \nabla \phi_k.
\]
Here, we use the fact that $U_i$ and $V_s$ are smooth functions and 
take their values at some points within $\omega$.
Similarly,
\begin{equation}
\begin{split}
\int_{R_\omega} \kappa \nabla (\phi_i^m \nabla_m U_i)  \cdot \nabla v\approx
\nabla_m U_i (x_\omega)\int_{R_\omega^+} \kappa \nabla \phi_i^m\cdot \nabla v =\\
\nabla_m U_i (x_\omega)  \nabla_k V_s (x_\omega) \int_{R_\omega} \kappa \nabla \phi_i^m \cdot \nabla \phi_s^k +
\nabla_m U_i (x_\omega)  V_s (x_\omega) \int_{R_\omega} \kappa \nabla \phi_i^m \cdot \nabla \phi_s\\
\nabla_m U_i (x_\omega)  \nabla_k V_s (x_\omega) \alpha_{is}^{km} + \nabla_m U_i (x_\omega)  V_s (x_\omega) \beta_{is}^{m*},\\
\end{split}
\end{equation}
where
\[
\alpha_{is}^{km}= \int_{R_\omega} \kappa \nabla \phi_i^m \cdot \nabla \phi_s^k.
\]
Next, using continuous approximations for $U_i$ and $V_i$, we can write
\begin{equation}
\begin{split}
\int_{R_\omega} \kappa \nabla u  \cdot \nabla v \approx
U_i  \beta_{ij}^{n*} \nabla_n V_j 
+
 U_i  \beta_{ij}^*   V_j +\\
  \nabla_m U_i\alpha_{ij}^{mn} \nabla_n V_j
+ 
\nabla_m U_i \beta_{ij}^{m*} V_j.
\end{split}
\end{equation}

Note that the definitions of $\alpha$'s and $\beta$'s
are using the volume of $R_\omega$ (which is of 
the same order as $R_\omega^+$). Moreover,
we also have the following scalings.
Assume $\epsilon$ is a diameter of RVE.
First, we note that
\begin{equation}
\label{eq:scalings}
\begin{split}
\|\phi_i\|=O(1),\ \|\nabla \phi_i\|=O({1\over\epsilon})\\
\|\phi_i^m\|=O(\epsilon),\ \|\nabla \phi_i^m\|=O(1).\\
\end{split}
\end{equation}
Using the formulas for $\alpha$'s and $\beta$'s, we
have the following scalings. 
\begin{align*}
\beta_{ij}^{m*}  =O(\cfrac{|R_\omega|}{\epsilon}),\ \alpha_{ij}^{mn}=O(|R_\omega|),\ 
\beta_{ij}^*  =O(\cfrac{|R_\omega|}{\epsilon^{2}}).
\end{align*}

We then define  rescaled quantities
 $\widehat{\beta}_{ij}$, 
$\widehat{\alpha}_{ij}$,
$\widehat{\beta}_{ij}$, 
$\widehat{\alpha}_{ij}$ such that
\begin{equation}
\begin{split}
\widehat{\beta}_{ij}  =\cfrac{|R_\omega|}{\epsilon^{2}}\beta_{ij}^*,\ 
\widehat{\beta}_{ij}^{n*}  =\cfrac{|R_\omega|}{\epsilon}\beta_{ij}^{n},\ 
\widehat{\alpha}_{ij}^{mn}  =|R_\omega|\alpha_{ij}^{mn}.
\end{split}
\end{equation}

With these scaling, we have
\begin{equation}
\label{eq:mc1}
\begin{split}
\int_\Omega \kappa \nabla u \cdot \nabla v \approx
 \int_\Omega \widehat{\alpha_{ij}^{mn}} \nabla_m U_i \nabla_n V_j + \\
{1\over\epsilon}\int_\Omega \widehat{\beta_{ij}^m} \nabla_m U_i V_j +
{1\over\epsilon}\int_\Omega  \widehat{\beta_{ij}^m} U_i \nabla_m V_j +
{1\over\epsilon^2}   \int_\Omega \widehat{\beta_{ij}}  U_i V_j.
\end{split}
\end{equation}
The sum of the second and third terms is negligible (this
can be shown by integration by parts).
It can be shown that
\[
\sum_j \widehat{\beta_{ij}}\approx 0.
\]
The last term can be written as
\begin{equation}
\begin{split}
 \int_\Omega \widehat{\beta_{ij}}  U_i V_j= \sum_{j\not = i} \int_\Omega \widehat{\beta_{ij}}  (U_i - U_j)V_j,
\end{split}
\end{equation}
which gives a form that is often used in multicontinuum 
models.
If we ignore the second and the third term in (\ref{eq:mc1}), we get
\begin{equation}
\begin{split}
-\nabla_n (\widehat{\alpha_{ij}^{mn}} \nabla_m U_j) +
   {1\over \epsilon^2}\widehat{\beta_{ij}} U_j=f_i
\end{split}
\end{equation}

We see from the equation that the reaction term
is dominant unless we deal with large diffusions 
(high contrast). When reaction terms dominate,
we have all $U_i$'s are approximately the same. 
 Thus, in general, to define
appropriate multicontinuum models (when $U_i$'s differ), one needs 
appropriate multicontinuum definitions, which
we will do in Section \ref{sec:spectral}.

If we have one continuum (as in standard homogenization),
then $\phi_1=1$ and $\beta_{ij}^p=0$.
The function $\phi_i^m$ will have the averages
$\int_{R_\omega^p} \phi_1^m=\int_{R_\omega^p} (x_m - x_m^0)$,
where $x_m^0=\langle x_m \rangle_{R_\omega^{p_0}}$. In this 
regard, $\phi_1^m$ acts as having linear growth as
in homogenization case (see later discussions
on imposing gradient constraints).

\subsection{Time dependent case}

We consider
\[
\int_{t_n}^{t_{n+1}} \int_\Omega u_t v + \int_{t_n}^{t_{n+1}} \int_\Omega \kappa \nabla u \cdot \nabla v = \int_{t_n}^{t_{n+1}} \int_\Omega f v.
\]
We introduce cell problems
(we keep the same notations as the stationary case)
 as time-dependent cell problems in 
$R_\omega^+\times [t_n, t_{n+1}]$
\begin{equation}
\label{eq:cell_grad_time}
\begin{split}
\int_{t_n}^{t_{n+1}} \int_{R_\omega^+} (\phi^m_i)_t v +\int_{t_n}^{t_{n+1}} \int_{R_\omega^+}\kappa \nabla \phi^m_i\cdot \nabla v -  
\\
\sum_p {\beta_{ij}^{mp}\over \int_{t_n}^{t_{n+1}}\int_{R_\omega}\psi_j} \int_{t_n}^{t_{n+1}}  
\int_{R_\omega}\psi_j  v =0\\
\int_{t_n}^{t_{n+1}}  \int_{R_\omega^p}  \phi^m_i \psi_j = \delta_{ij}  \int_{t_n}^{t_{n+1}}  \int_{R_\omega^p}  (x_m-c_{mj}) \psi_j\\
\int_{t_n}^{t_{n+1}}  \int_{R_\omega^{p_0}}  (x_m-c_{mj}) \psi_j=0\\
\phi_i^m(t=t_n)=\xi_i^m.
\end{split}
\end{equation}
The second cell problem is defined with constraints on the average
solutions in $R_\omega^+\times [t_n, t_{n+1}]$
\begin{equation}
\label{eq:cell_aver_time}
\begin{split}
\int_{t_n}^{t_{n+1}} \int_{R_\omega^+} (\phi_i)_t v +
\int_{t_n}^{t_{n+1}} \int_{R_\omega^+}\kappa \nabla \phi_i\cdot \nabla v -  \\
\sum_p {\beta_{ij}^p\over \int_{t_n}^{t_{n+1}}\int_{R_\omega}\psi_j} \int_{t_n}^{t_{n+1}} \int_{R_\omega}\psi_j  v =0\\
\int_{t_n}^{t_{n+1}} \int_{R_\omega^p}  \phi_i \psi_j = \delta_{ij} \int_{t_n}^{t_{n+1}} \int_{R_\omega^p}   \psi_j \\
\phi_i(t=t_n)=\psi_i.
\end{split}
\end{equation}

We again assume that in $R_\omega=R_\omega^{p_0}$
\begin{equation}
\begin{split}
u\approx  \phi_i U_i  + \phi_i^m \nabla_m U_i.
\end{split}
\end{equation}
Next, we derive multicontinuum equations for $U_i$.
Then, we have (for any $v=\phi_s V_s  + \phi_s^k \nabla_k V_s$)
\begin{equation}
\label{eq:assump2}
\begin{split}
\int_{t_n}^{t_{n+1}} \int_\Omega f v =
\int_{t_n}^{t_{n+1}}\int_\Omega u_t v +
 \int_{t_n}^{t_{n+1}}\int_\Omega \kappa \nabla u \cdot \nabla v = \\
\sum_\omega \int_{t_n}^{t_{n+1}} \int_\omega u_t v +
\sum_\omega \int_{t_n}^{t_{n+1}} \int_\omega \kappa \nabla u \cdot \nabla v \approx\\
\sum_\omega {|\omega|\over |R_\omega|} \int_{t_n}^{t_{n+1}} \int_{R_\omega}  u_t v +  \sum_\omega {|\omega|\over |R_\omega|} \int_{t_n}^{t_{n+1}} \int_{R_\omega} \kappa \nabla u \cdot \nabla v.
\end{split}
\end{equation}
Next, 
\begin{equation}
\label{eq:exp1}
\begin{split}
\int_{t_n}^{t_{n+1}}\int_{R_\omega} u_t v + \int_{t_n}^{t_{n+1}}\int_{R_\omega} \kappa \nabla u  \cdot \nabla v \approx
(U_i)_t V_k \int\int \phi_i \phi_k + \\
U_i V_k \int\int \left((\phi_i)_t \phi_k + \kappa \nabla \phi_i \cdot \nabla \phi_k \right) +
\nabla_m U_i \nabla_n V_k \int\int \left((\phi^m_i)_t \phi^n_k + \kappa \nabla \phi_i^m \cdot \nabla \phi_k^n \right) +\\
U_i \nabla_n V_k \int\int \left((\phi_i)_t \phi^n_k + \kappa \nabla \phi_i \cdot \nabla \phi_k^n \right) +
\nabla_m U_i  V_k \int\int \left((\phi_i^m)_t \phi_k + \kappa \nabla \phi_i^m \cdot \nabla \phi_k\right) +\\
(U_i)_t\nabla_n V_k \int\int \phi_i\phi_k^n= m_{ik} (U_i)_t V_k+\beta_{ik}U_i V_k +\alpha_{ij}^{nm} \nabla_m U_i \nabla_n V_k + \\
\beta_{ik}^{n} U_i \nabla_n V_k + \alpha_{ik}^{m} \nabla_m U_i  V_k +m_{ij}^n (U_i)_t\nabla_n V_k.
\end{split}
\end{equation}
We note that if the cell problems do not contain $t$, i.e.,
$\phi_i^m$ and $\phi_i$ do not depend on $t$,
we get similar homogenized equations (which can easily be derived
from (\ref{eq:exp1})).
Here for simplicity of the notations, $\int\int \cdot =  \int_{t_n}^{t_{n+1}} \int_{R_\omega}\cdot$ and the notations for $m_{ik}$, $\beta_{ik}$,
$\alpha_{ij}^{nm}$, $\beta_{ik}^{n}$, $\alpha_{ik}^{m}$, and $m_{ij}^n$
can be seen from the above equality.
 We neglect $\phi_i^m \nabla_m U_i$ and
$(U_i)_t \nabla_n V_k \int\int \phi_i \phi_k^n$
based on scalings (\ref{eq:scalings1}).

Assume $\epsilon$ is a diameter of RVE.
We note that
\begin{equation}
\label{eq:scalings1}
\begin{split}
\|\phi_i\|=O(1),\ \|\nabla \phi_i\|=O({1\over\epsilon})\\
\|\phi_i^m\|=O(\epsilon),\ \|\nabla \phi_i^m\|=O(1).\\
\end{split}
\end{equation}
Thus, the last term in (\ref{eq:exp1}) can be neglected.

Using continuous approximation $U_i$'s and $V_i$'s,
we can write
\begin{equation}
\begin{split}
\int_{t_n}^{t_{n+1}} \int_\Omega  f v=
\int_{t_n}^{t_{n+1}} \int_\Omega u_t v + \int_{t_n}^{t_{n+1}} \int_\Omega \kappa  \nabla u  \cdot \nabla v  \approx\\
\int_{t_n}^{t_{n+1}}  \int_\Omega \widetilde{m_{ij}} (U_i)_t V_j+
\int_{t_n}^{t_{n+1}}  \int_\Omega \widetilde{\alpha_{ij}^{mn}} \nabla_m U_i \nabla_n V_j + \\
\int_{t_n}^{t_{n+1}}  \int_\Omega \widetilde{\alpha_{ij}^m} \nabla_mU_i V_j +
\int_{t_n}^{t_{n+1}} \int_\Omega  \widetilde{\beta_{ij}^n} U_i \nabla_n V_j +\\
\int_{t_n}^{t_{n+1}}  \int_\Omega  \widetilde{\beta_{ij}} U_i V_j,
\end{split}
\end{equation}
where $\widetilde{\cdot}$ denotes rescaled $\cdot$ with scaling of $R_\omega$
and $\Delta t$ (i.e., $\widetilde{\cdot} = { \cdot \over |R_\omega||\Delta t|}$
so that we can write the integrals.

Formally, we can write the system of differential equations
\begin{equation}
\begin{split}
(U_i)_t - \nabla_n \widetilde{\alpha_{ij}^{mn}} \nabla_m U_j +
 \widetilde{\beta_{ij}^m} \nabla_mU_j  \\
- \nabla_n(\widetilde{\alpha_{ij}^n} U_j)  +
   \widetilde{\beta_{ij}} U_j=f_i
\end{split}
\end{equation}

Using the formulas for $\alpha$'s and $\beta$'s, we
have the following scalings (if we ignore the terms with temporal derivatives)
\begin{align*}
\beta_{ij}^{m} & =O(\cfrac{|R_\omega|\Delta t}{\epsilon}),\alpha_{ij}^{mn}=O(|R_\omega|\Delta t)\\
\beta_{ij} & =O(\cfrac{|R_\omega|\Delta t}{\epsilon^{2}}),\alpha_{ij}^{n}=O(\cfrac{|R_\omega|\Delta t}{\epsilon}).
\end{align*}
One can make similar argument as in steady state case.

\section{Choices of continua. Spectral continua}
\label{sec:spectral}

In this section, we discuss how high contrast can balance
$\alpha$'s and $\beta$'s. 
First, we note that if reaction terms (represented via $\beta$)
dominate, then all $U_i$'s are approximately the same and we do not
have multicontinuum (i.e., different averages in different continua).
We assume steady-state case and
two continua, where the continuum $1$ has high-contrast $\kappa=O(\eta)$,
$\eta$ is large, and the continum $2$ has a conductivity of order $1$.
The next arguments do not take into account RVE sizes and are purely in terms
of $\eta$.
We can see that $\alpha_{is}^{km}$ (the diffusivity) 
is large $O(\eta)$ (at least in some direction)
since the local solutions have linear growth conditions. In general, 
the scalings of $\alpha$'s in terms of the contrast depend
on heterogeneities (see numerical results).
On the other hand, $\beta_{11}$ is of order $1$ (in terms of the contrast, while it depends on the size of $R_\omega$).
Since $\sum_{j}\beta_{ij} \approx 0$, we can conclude that other
$\beta$'s are of order $1$. Thus, if the contrast balances the RVE size
(e.g., $\eta= O(\epsilon^{-2}))$ (where $\epsilon$ is the size of $R_\omega$), 
we expect that $\beta$ terms do not dominate
and there are differences between average states and one has multicontinuum
homogenized limit.

Next, we discuss how one can identify the continuum via local spectral
problems. 
We consider $\beta_{ij}=\int_{R_\omega} \kappa \nabla \phi_i \cdot \nabla \phi_j$ and  assume, for simplicity, that the cell problem (\ref{eq:cell_aver_aver}) is 
formulated in $R_\omega$.
We would like to minimize $\beta_{ij}$ with constraints.
We assume (1) $\phi_i$'s are in the space of local solutions ($\zeta_j$)
(2) $\psi_j$ can take any values and the constraints are given by
$\int_{R_\omega} \kappa \phi_i \psi_j = \delta_{ij}$.
The functions $\zeta_j$ are local homogeneous solutions,
$\nabla \cdot (\kappa \nabla \zeta_j)=0$ in $R_\omega$, with boundary conditions
$\zeta_j = \delta^h_j(x)$ on $\partial R_\omega$, where $\delta_j(x)$ is a fine-grid hat function defined on the boundary of $R_\omega$.
In this case, if we consider the eigenvalue problem 
\begin{equation}
\label{eq:spectral1}
-div(\kappa \nabla \eta_j) = \lambda_j \kappa \eta_j.
\end{equation}
(see \cite{GMsFEM13,chung2018constraint}). 
with corresponding Rayleigh quotient,
\[
{\int_{R_\omega} \kappa |\nabla \phi|^2 \over \int_{R_\omega} \kappa |\phi|^2},
\]
it is clear that $\psi_j=\phi_j=\eta_j$.
The eigenvectors corresponding to the smallest eigenvalues are
constant functions in high-contrast regions.
In general, one can identify high-contrast regions by finding nearly
constant gradient regions of $\nabla \phi_j$.
Moreover, the number of smallest eigenvalues will correspond to
the number of high-conductivity channels (channels that connect boundaries of
$R_\omega$).

Based on the above discussion, one can perform local spectral decomposition
based on (\ref{eq:spectral1}) and identify $\psi_j$ based on smallest
eigenvalues
that correspond to high contrast (they scale as the inverse of high contrast).
Using these eigenvectors, the local problems (\ref{eq:cell_aver_aver})
and (\ref{eq:cell_grad_aver}) are solved. We can use instead of 
$\beta_{ij}^{mp}/\int_{R_\omega^p} \psi_j^p$ and
$\beta_{ij}^{p}/\int_{R_\omega^p} \psi_j^p$, the terms without denominator,
$\beta_{ij}^{mp}$ and  $\beta_{ij}^{p}$ so that not to worry that
$\psi_j^p$ may vanish.

\section{Remarks}
\subsection{Multicontinuum derivation based on average and gradient constraint problems in representative volumes.}

One can also use gradient type constraints in addition.
We demonstrate this and point out some issues in this procedure.
We consider a steady state diffusion equation.

In each RVE, we introduce 
two sets of cell problems with constraints
formulated in  $R_\omega$ (though, one can use oversampling).
\begin{equation}
\label{eq:cell_grad}
\begin{split}
\int_{R_\omega}\kappa \nabla \phi^m_i \cdot\nabla v -  {\alpha_{ij}^{mn}\over \int_{R_\omega}\psi_j} \int_{R_\omega}\psi_j \nabla_n v -  
{\beta_{ij}^{m}\over \int_{R_\omega}\psi_j} \int_{R_\omega}\psi_j  v =0\\
\int_{R_\omega} \nabla \phi^m_i \psi_j = \delta_{ij} e_m \int_{R_\omega} \psi_j\\
\int_{R_\omega} \phi^m_i \psi_j = 0,\\
\end{split}
\end{equation}
where $e_m$ is $m$th unit vector.
The second cell problem is defined with constraints on the average
solutions in $R_\omega$ 
\begin{equation}
\label{eq:cell_aver}
\begin{split}
\int_{R_\omega}\kappa \nabla \phi_i\cdot \nabla v -  {\alpha_{ij}^{n}\over \int_{R_\omega}\psi_j} \int_{R_\omega}\psi_j \nabla_n v - 
{\beta_{ij}\over \int_{R_\omega}\psi_j} \int_{R_\omega}\psi_j  v =0\\
\int_{R_\omega}  \phi_i \psi_j = \delta_{ij}  \int_{R_\omega} \psi_j \\
\int_{R_\omega}  \nabla \phi_i \psi_j = 0. \\
\end{split}
\end{equation}

Some properties of Lagrange multipliers are discussed.
We note 
$\sum_j \beta_{ij}=0$,
which can be obtained by taking $v=1$ in (\ref{eq:cell_aver}).
If we take $v=\phi_j$ in  (\ref{eq:cell_grad}) then
$\int_{R_\omega} \kappa \nabla \phi_i^m \cdot \nabla \phi_j = \beta_{ij}^m$.
If we take $v=\phi_k^m$ in  (\ref{eq:cell_aver}), we get
$\int_{R_\omega} \kappa \nabla \phi_i \cdot \nabla \phi_k^m  = \alpha_{ik}^m.$
Therefore,
$\alpha_{ij}^m= \beta_{ij}^m$.
If we take $v=\phi_k^s$ in (\ref{eq:cell_grad}), then,
we have
$\int_{R_\omega} \kappa \nabla \phi_i^m  \cdot \nabla\phi_k^s = \alpha_{ik}^{ms}.$
If
 we take $v=\phi_s$ in (\ref{eq:cell_aver}),
we get 
$ \int_{R_\omega} \kappa \nabla \phi_i \cdot \nabla \phi_s = \beta_{is}.$

Next, we assume that the local solution in $\omega$ can be represented
by constraint problems in  $R_\omega$
\begin{equation}
\begin{split}
u\approx  \phi_i U_i  + \phi_i^m V_i^m,
\end{split}
\end{equation}
where
\begin{equation}
\begin{split}
U_i ={\int_{R_\omega} u \psi_i\over \int_{R_\omega}\psi_i},\ \ \  
V_i^m={\int_{R_\omega} \nabla_m u \psi_i \over  \int_{R_\omega}\psi_i}.
\end{split}
\end{equation}
Similarly, we introduce for test functions
\begin{equation}
\begin{split}
P_i ={\int_{R_\omega} v \psi_i\over \int_{R_\omega}\psi_i},\ \ \  
Q_i^m={\int_{R_\omega} \nabla_m v \psi_i \over  \int_{R_\omega}\psi_i}.
\end{split}
\end{equation}
For any $v\in H_0^1$, we have
\begin{equation}
\label{eq:assump3}
\begin{split}
\int_\Omega f v = \int_\Omega \kappa \nabla u  \cdot \nabla v = \sum_\omega \int_\omega \kappa \nabla u  \cdot \nabla v \approx
\sum_\omega {|\omega|\over |R_\omega|} \int_{R_\omega} \kappa \nabla u  \cdot \nabla v, 
\end{split}
\end{equation}
where we make an assumption that integrated average over
RVE can represent the whole 
computational element $\omega$.
Next, we approximate each term 
\begin{equation}
\begin{split}
\int_{R_\omega} \kappa \nabla u   \cdot \nabla v =
\int_{R_\omega} \kappa \nabla (\phi_i U_i)  \cdot \nabla v  + 
\int_{R_\omega} \kappa \nabla (\phi_i^m V_i^m)  \cdot \nabla v. 
\end{split}
\end{equation}
We assume that the variation of $U_i$ and $V_i^m$
are small compared to
the variations of $\phi_i$ and $\phi_i^m$ (since
they vary at RVE scale)
and assume
$\int_{R_\omega} \kappa \nabla (\phi_i U_i)  \cdot \nabla v\approx
\int_{R_\omega} \kappa (\nabla \phi_i) U_i  \cdot \nabla v$
and
$\int_{R_\omega} \kappa \nabla (\phi_i^m V_i^m )   \cdot \nabla v
\approx
\int_{R_\omega} \kappa \nabla (\phi_i^m) V_i^m  \cdot \nabla v$. 
Then (denoting $x_\omega$ center of $R_\omega$),
\begin{equation}
\begin{split}
\int_{R_\omega} \kappa \nabla (\phi_i U_i)  \cdot \nabla v \approx
U_i(x_{\omega}) \int_{R_\omega} \kappa \nabla \phi_i  \cdot \nabla v =\\
U_i(x_{\omega}) 
 {\alpha_{ij}^{n}\over \int_{R_\omega}\psi_j} \int_{R_\omega}\psi_j  \cdot \nabla_n v 
+
 U_i(x_{\omega})   {\beta_{ij}\over \int_{R_\omega}\psi_j} \int_{R_\omega}\psi_j  v=\\
U_i(x_{\omega}) \alpha_{ij}^{n} Q_j^n (x_\omega)
+
 U_i(x_{\omega})  \beta_{ij}   P_j (x_\omega).
\end{split}
\end{equation}
Similarly,
\begin{equation}
\begin{split}
\int_{R_\omega} \kappa \nabla (\phi_i^m \nabla_m U_i)  \cdot  \nabla v\approx
 V_i^m (x_\omega)\int_{R_\omega} \kappa \nabla \phi_i^m  \cdot \nabla v =\\
 {\alpha_{ij}^{mn}\over \int_{R_\omega}\psi_j} V_i^m (x_\omega) \int_{R_\omega}\psi_j \nabla_n v + 
{\beta_{ij}^{m}\over \int_{R_\omega}\psi_j} V_i^m (x_\omega) \int_{R_\omega}\psi_j  v =\\
  V_i^m (x_\omega) \alpha_{ij}^{mn} Q_j^n(x_\omega) 
+ 
 V_i^m (x_\omega) \beta_{ij}^{m} P_j(x_\omega).
\end{split}
\end{equation}
Next, using continuous approximations for all quantities, we have
\begin{equation}
\begin{split}
\int_{R_\omega} \kappa \nabla u  \cdot \nabla v \approx
U_i  \alpha_{ij}^{n} Q_j^n 
+
 U_i  \beta_{ij}   P_j +\\
  V_i^m \alpha_{ij}^{mn} Q_j^n
+ 
V_i^m \beta_{ij}^{m} P_j.
\end{split}
\end{equation}
Consequently, 
\begin{equation}
\label{eq:multi}
\begin{split}
\int_\Omega f v =\int_\Omega \kappa \nabla u  \nabla v \approx
 \int_\Omega \widetilde{\alpha_{ij}^{mn}} V_i^m Q_j^n + \\
\int_\Omega \widetilde{\beta_{ij}^m} V_i^m P_j +
\int_\Omega  \widetilde{\alpha_{ij}^n} U_i Q_j^n +
 \int_\Omega  \widetilde{\beta_{ij}} U_i P_j,
\end{split}
\end{equation}
where $\widetilde{\cdot}$ denotes rescaled $\cdot$ with scaling 
(i.e., $\widetilde{\cdot} = { \cdot \over |R_\omega|}$).
The equation (\ref{eq:multi}) is a multicontinuum model
equation. 
Similar scaling arguments as before can be applied.

The resulting equations are similar to those we obtained earlier. 
In our numerical studies, we have found that 
the cell problems using multiple simultaneous constraints (e.g., average 
solution and average gradient) are prone to large errors. This is 
because one needs to choose these simultaneous constraints
(for averages and gradients) such that the local problems
$\phi_i$ and $\phi_i^m$ have similar fine-scale features as the exact solution.
When there is no compatibility (i.e., averages and gradients do not correspond to each other), then there are large errors,
especially at the interfaces of the continuum, which can cause
large errors on average characteristics. For this reason, 
 we have found one constraint
cell problems to be more accurate, easy to implement, and
easy to remove boundary effects with oversampling.

\subsection{Nonlinear case. Steady state}

The derivations presented earlier can be done for nonlinear
problems. 
We consider
\[
\int_\Omega \kappa(x, \nabla u) \cdot \nabla v= \int_\Omega f v.
\]
In this case, we can split the cell problem into average-based
and gradient based due to nonlinear interaction and we consider
the following cell problem.
\begin{equation}
\begin{split}
\int_{R_{\omega}^{+}}\kappa(x,\nabla\phi(\eta,\xi))\cdot \nabla v-\sum_{p}\alpha_{\omega}^{p,j}(\eta,\xi)\int_{R_{\omega}^{p}}\psi_{j}^{p}v  =0\\
\int_{R_{\omega}^{p}}\phi^{\eta,\xi}\psi_{j}^{p} =\eta_{j}\int_{R_{\omega}^{p}}\psi_{j}^{p}+\xi_{j}\cdot\int_{R_{\omega}^{p}}(x-c)\psi_{j}^{p}.
\end{split}
\end{equation}
We will use the following approximations 
$$u\approx\phi(U,\nabla U)$$
$$v|_{R_{\omega}^{p}}\approx\sum_{i}V_{i}\cfrac{\psi_{i}^{p}}{\int\psi_{i}^{p}}+\sum_{i}\cfrac{\psi_{i}^{p}}{\int\psi_{i}^{p}}(x-c_{i}) \cdot \nabla V_{i}.$$
Then,
\begin{equation}
\begin{split}
\int_{R_\omega^+}\kappa(x,\nabla\phi(U,\nabla U))\cdot \nabla v  =\sum_{p}\alpha_{\omega}^{p,j}(U,\nabla U)V_{i}\cfrac{\int_{R_{\omega}^{p}}\psi_{j}^{p}\psi_{i}^{p}}{\int_{R_{\omega}^{p}}|\psi_{i}^{p}|^{2}}+\\
\sum_{p}\alpha_{\omega}^{p,j}(U,\nabla U)\cfrac{\int_{R_{\omega}^{p}}\psi_{j}^{p}\psi_{i}^{p}(x-c_{i})}{\int_{R_{\omega}^{p}}\psi_{i}^{p}}\cdot \nabla V_{i}
  =\\
\sum_{p}\alpha_{\omega}^{p,j}(U,\nabla U)V_{j}+\sum_{m,p}\beta^{jm}(U,\nabla U)\nabla_{m}V_{j},
\end{split}
\end{equation}
where all $U$'s and $V$'s are taken to be constants at RVE-level and
\begin{align*}
\gamma^{j}(\eta,\xi) & =\sum_{p}\alpha_{\omega}^{p,j}(\eta,\xi)\\
\beta^{jm}(\eta,\xi) & =\sum_{p}\alpha_{\omega}^{p,j}(\eta,\xi)\cfrac{\int_{R_{\omega}^{p}}(x_{m}-c_{m})\psi_{j}^{p}}{\int\psi_{j}^{p}}.
\end{align*}
Here, we took account
$
\int_{R_{\omega}^{p}}\psi_{j}^{p}\psi_{i}^{p}=\delta_{ij}\int_{R_{\omega}^{p}}\psi_{i}^{p}$ and $\int_{R_{\omega}^{p}}\psi_{j}^{p}\psi_{i}^{p}(x-c_{i})=\delta_{ij}\int_{R_{\omega}^{p}}\psi_{i}^{p}(x-c_{i})$.
We get the following multicontinuum equations
\[
\int_\Omega\gamma^{j}(U,\nabla U)V_{j}+\int_\Omega\beta^{j}(U,\nabla U)\cdot\nabla V_{j}=f_{j}.
\]
In general, these equations are complicated to solve. 
Some machine learning techniques are needed to train the
local upscaled quantities.
We have presented some cases in \cite{leung2019space,vasilyeva2020learning,chung2018nonlinear}.
In some
special cases, one can simplify the resulting multicontinuum equations.

\section{Numerical results}

In this section, we will present numerical examples of the proposed
upscaling method. We will present three numerical examples. The goal
is to show that our proposed algorithm is accurate and the cell
problem solutions provide better accuracy as we increase oversampling size.

\begin{figure}
\centering
\includegraphics[scale=0.3]{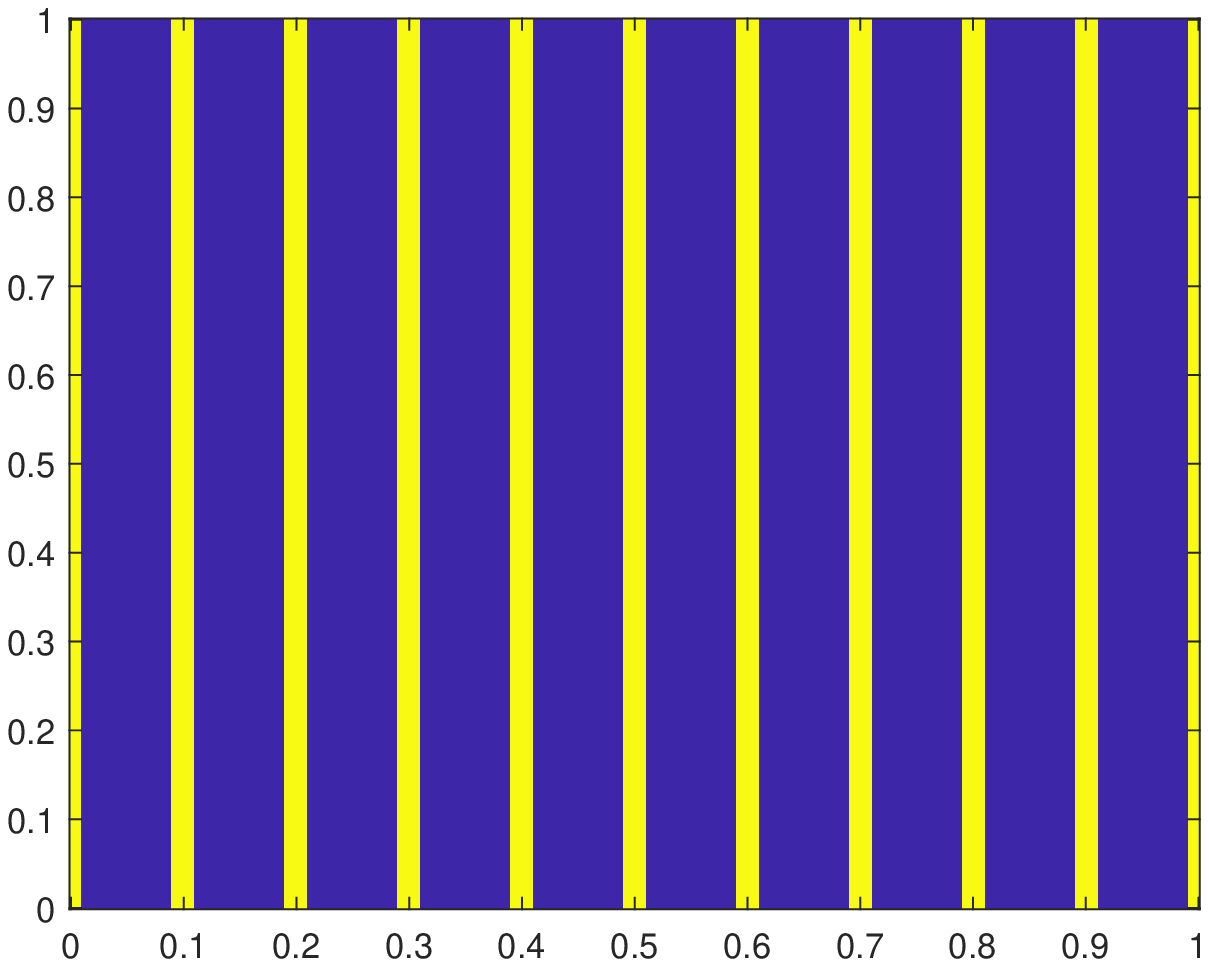}
\includegraphics[scale=0.3]{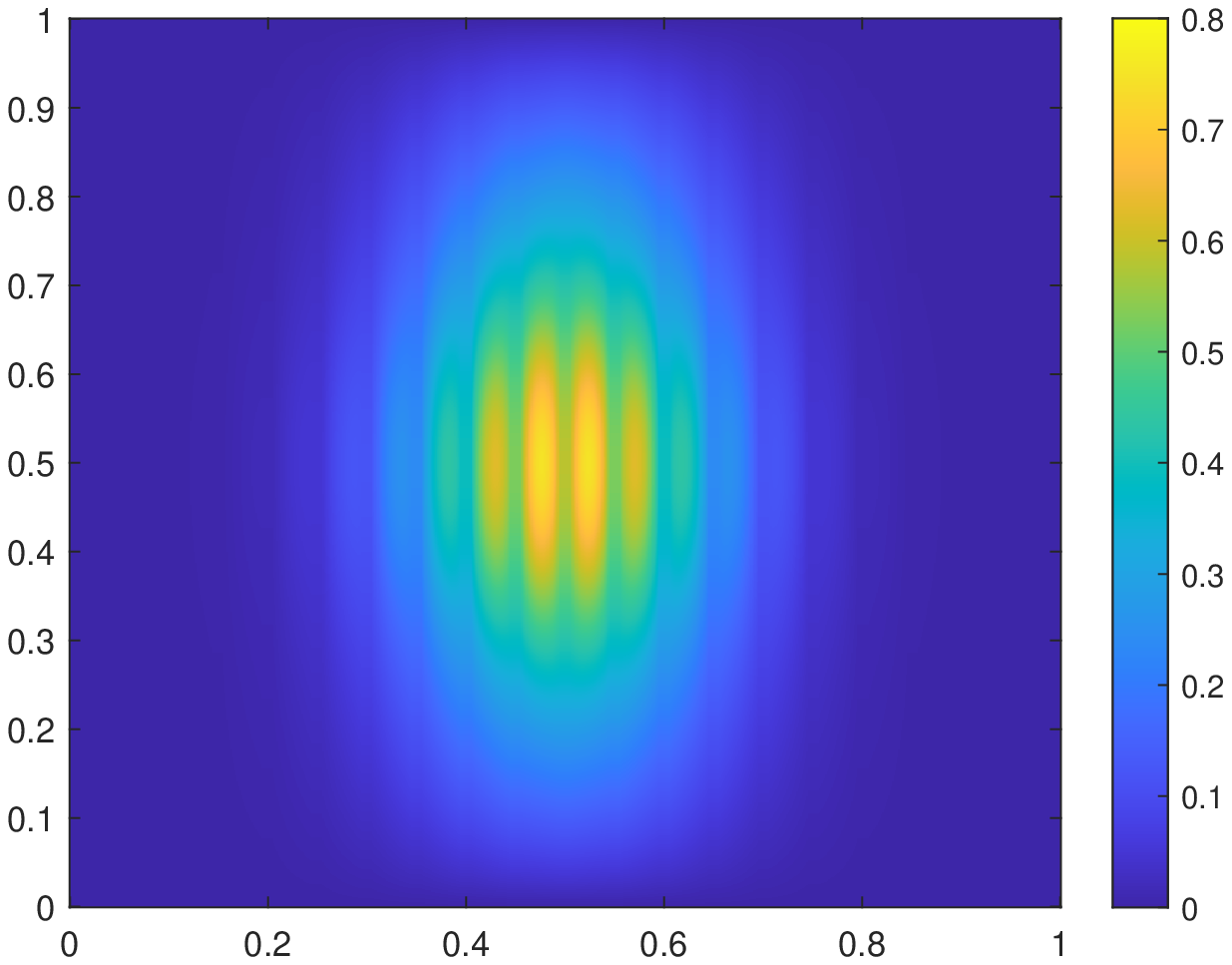}
\caption{Left: $\kappa$ for Case 1  with $\epsilon=\cfrac{1}{10}$. Right: The solution snapshot.}
\label{fig:perm1}
\end{figure}

In the first example, we consider the layered medium parameter $\kappa$
(see Figure \ref{fig:perm1}).
The period of $\kappa$ is denoted as $\epsilon$. We denote the low
conductivity region and the high conductivity region of $\kappa$
by $\Omega_{1}$ and $\Omega_{2}$, respectively.
The source term
$f$ and conductivity $\kappa$ is given as
\[
f(x)=\begin{cases}
1000\min\{\kappa\}e^{-40|(x-0.5)^{2}+(y-0.5)^{2}|} & x\in\Omega_{1}\\
e^{-40|(x-0.5)^{2}+(y-0.5)^{2}|} & x\in\Omega_{2}
\end{cases}
\]
and
\[
\kappa(x)=\begin{cases}
\cfrac{\epsilon}{10000} & x\in\Omega_{1}\\
\cfrac{1}{100\epsilon} & x\in\Omega_{2}
\end{cases}
\]
 The computational domain $\Omega$ is partitioned into $M\times M$
coarse grid. The coarse mesh size $H$ is defined as $H=1/M$. For
simplicity, we consider the whole coarse grid element as an RVE for
the corresponding coarse element. The oversampling RVE $R_\omega^{+}$
(or $\omega^+$)
for each coarse RVE $\omega$ is defined as an extension of $K$ by
$l$ layers of coarse grid element, where $l$ will be varied.

We will define the relative $L^{2}$- error in $\Omega_{1}$ and the
relative $L^{2}$- error in $\Omega_{2}$ by 
\[
e_{2}^{(i)}=\cfrac{\sum_{K}|\cfrac{1}{|K|}\int_{K}U_{i}-\cfrac{1}{|K\cap\Omega_{i}|}\int_{K\cap\Omega_{i}}u|^{2}}{\sum_{K}|\cfrac{1}{K\cap\Omega_{i}}\int_{K\cap\Omega_{i}}u|^{2}}.
\]
$K$ denotes the RVE, which is taken to be $\omega$.
This represents the $L_2$ error of our proposed approach.

Our goals in this section are the following.
\begin{itemize}

\item We show that the errors between the upscaled solutions and corresponding fine-grid solutions are small.

\item The errors are stable for different mesh sizes, RVE sizes, and contrasts. If the number of layers is appropriately chosen (to avoid boundary effects), 
the error will decrease as we decrease the mesh size.

\item We discuss the effective properties and show their values for different
values of mesh sizes and RVE sizes.

\item The numerical results show scalings of effective parameters.

\end{itemize}

\begin{figure}
\centering
\includegraphics[scale=0.3]{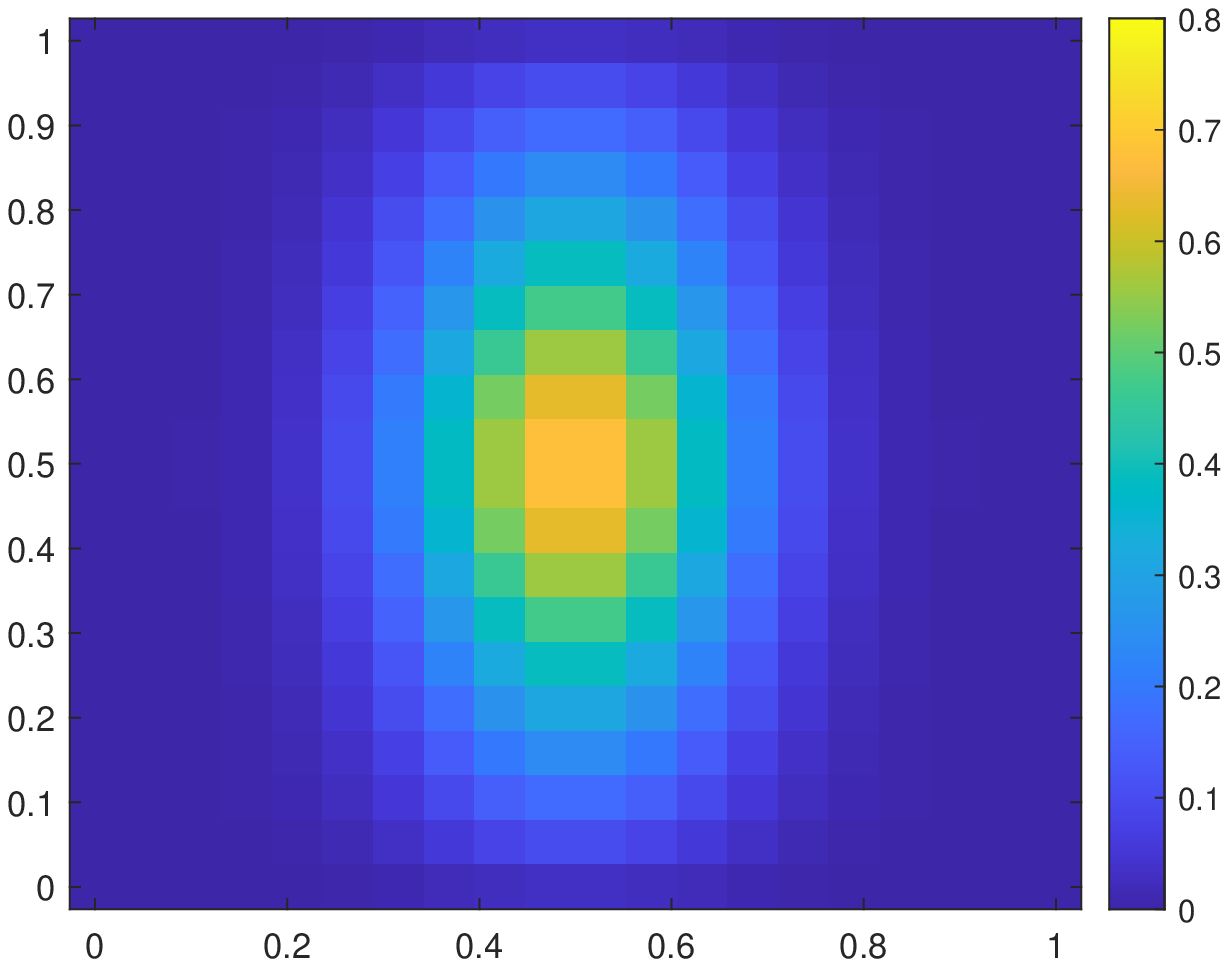} 
\includegraphics[scale=0.3]{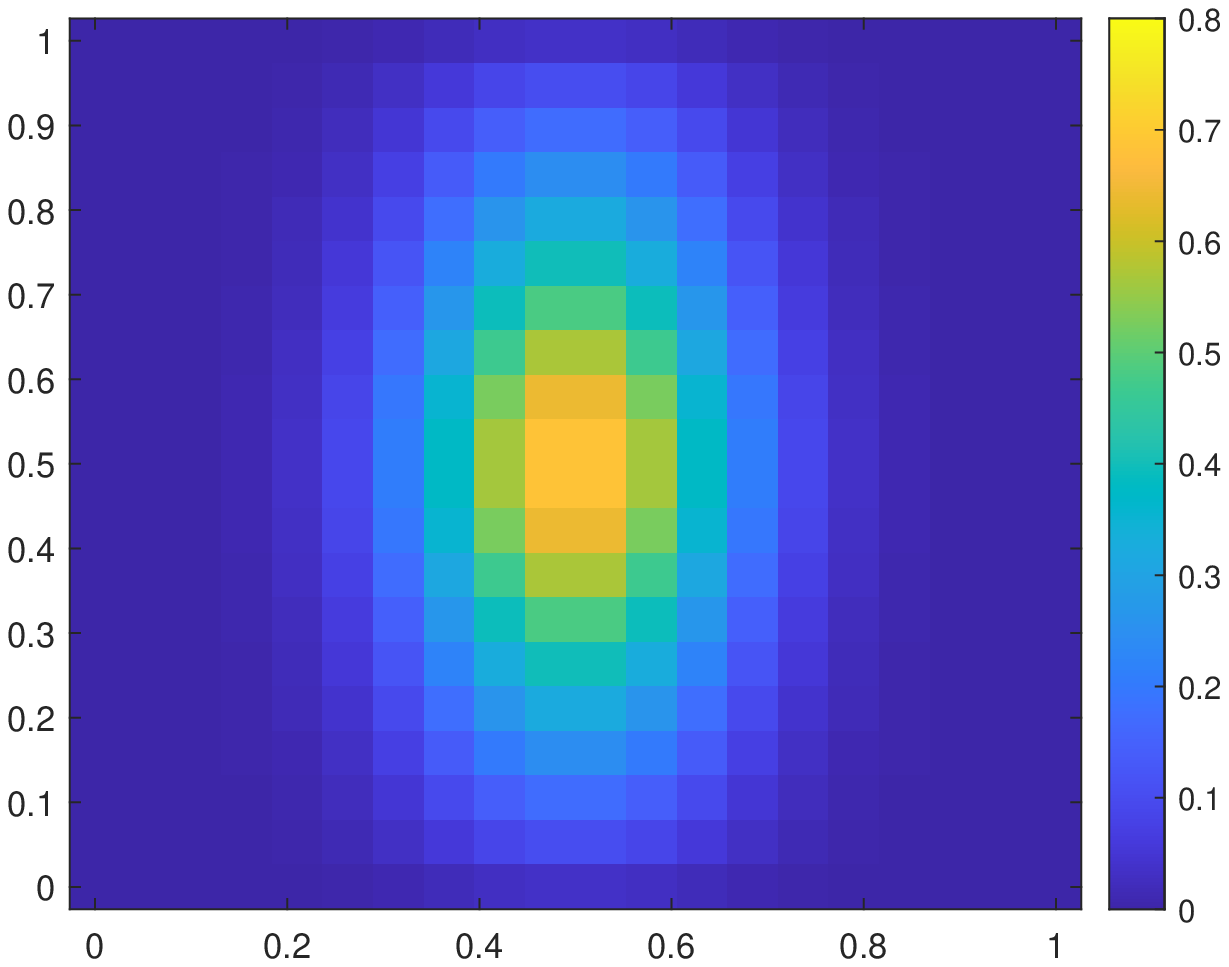}
\includegraphics[scale=0.3]{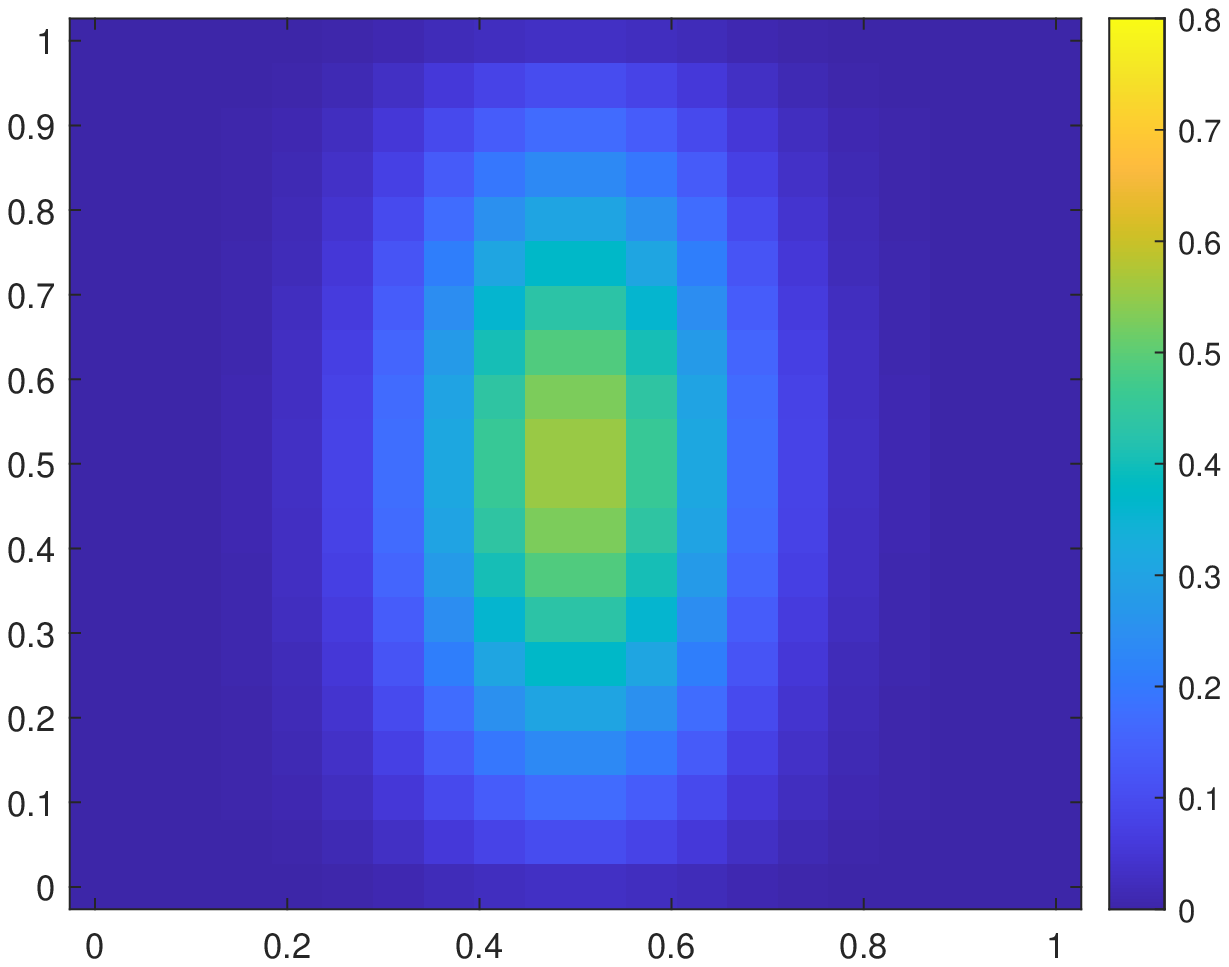} 
\includegraphics[scale=0.3]{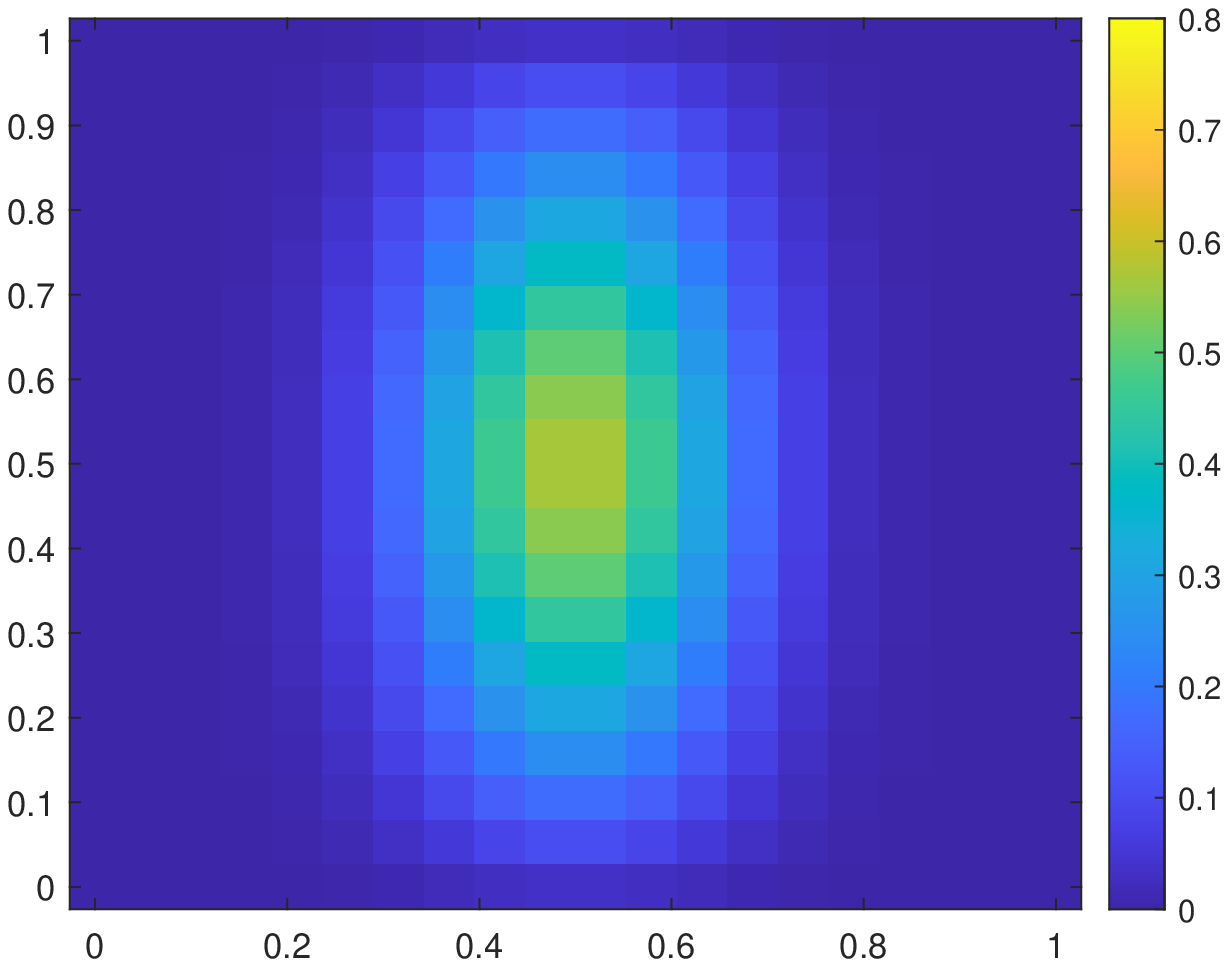}
\caption{Case 1. Top-Left:
Reference averaged solution in $\Omega_{1}$. Top-Right: Multiscale
solution $U_{1}$. Bottom-Left: Reference averaged solution in $\Omega_{2}$.
Bottom-Right: Multiscale solution $U_{2}$.}
\label{fig:compare_case1}
\end{figure}

For the first case, we present $e_2^{(i)}$ in Table \ref{tab:case1}.
We make several observations. First, we observe that the proposed
approach provides an accurate approximation of the averaged solution.
In Figure \ref{fig:compare_case1}, we depict upscaled solutions
and corresponding averaged fine-scale solutions. We observe that these
solutions are very close.
In the first table, we decrease the coarse-mesh size and it gets
closer to $\epsilon$. In standard numerical homogenization methods,
this was known to give a resonance error and error will increase.
Here, by choosing an appropriate number of layers, we observe that
the error remains small. In the second table, we decrease the period
size and observe that the error decreases to a certain level. In the third
table, we observe convergence as we decrease the mesh size and 
$\epsilon$. In general,
we expect a certain threshold error due to fine-scale discretization, which is
used to compute the solution.
In Table \ref{tab:case1_fixed}, we present the errors for fixed
contrast ratio $1/10000$ (in $\Omega_1$) and 1/10 (in $\Omega_2$).
As we decrease the mesh and RVE sizes, we observe that the upscaled
solution converges to the averaged fine-scale solution.

\begin{table}
\centering
\begin{tabular}{|c|c|c|c|c|}
\hline 
$l$ & $H$ & $\epsilon$ & $e_{2}^{(1)}$ & $e_{2}^{(2)}$\tabularnewline
\hline 
$\left\lceil -2\log(H)\right\rceil =5$ & $\cfrac{1}{10}$ & $\cfrac{1}{40}$ & 1.50\% & 1.37\%\tabularnewline
\hline 
$\left\lceil -2\log(H)\right\rceil =6$ & $\cfrac{1}{20}$ & $\cfrac{1}{40}$ & 0.48\% & 0.50\%\tabularnewline
\hline 
$\left\lceil -2\log(H)\right\rceil =8$ & $\cfrac{1}{40}$ & $\cfrac{1}{40}$ & 0.61\% & 0.60\%\tabularnewline
\hline 
\end{tabular} %
\begin{tabular}{|c|c|c|c|c|}
\hline 
$l$ & $H$ & $\epsilon$ & $e_{2}^{(1)}$ & $e_{2}^{(2)}$\tabularnewline
\hline 
$5$ & $\cfrac{1}{10}$ & $\cfrac{1}{10}$ & 4.60\% & 8.35\%\tabularnewline
\hline 
$5$ & $\cfrac{1}{10}$ & $\cfrac{1}{20}$ & 1.60\% & 1.31\%\tabularnewline
\hline 
5 & $\cfrac{1}{10}$ & $\cfrac{1}{40}$ & 1.50\% & 1.37\%\tabularnewline
\hline 
\end{tabular}

\begin{tabular}{|c|c|c|c|c|}
\hline 
$l$ & $H$ & $\epsilon$ & $e_{2}^{(1)}$ & $e_{2}^{(2)}$\tabularnewline
\hline 
$\left\lceil -2\log(H)\right\rceil =5$ & $\cfrac{1}{10}$ & $\cfrac{1}{10}$ & 4.60\% & 8.35\%\tabularnewline
\hline 
$\left\lceil -2\log(H)\right\rceil =6$ & $\cfrac{1}{20}$ & $\cfrac{1}{20}$ & 2.02\% & 2.40\%\tabularnewline
\hline 
$\left\lceil -2\log(H)\right\rceil =8$ & $\cfrac{1}{40}$ & $\cfrac{1}{40}$ & 0.61\% & 0.60\%\tabularnewline
\hline 
$\left\lceil -2\log(H)\right\rceil =9$ & $\cfrac{1}{80}$ & $\cfrac{1}{80}$ & 0.14\% & 0.14\%\tabularnewline
\hline 
\end{tabular}
\caption{Error comparison for Case 1.}
\label{tab:case1}
\end{table}

\begin{table}
\centering
\begin{tabular}{|c|c|c|c|c|}
\hline 
$l$ & $H$ & $\epsilon$ & $e_{2}^{(1)}$ & $e_{2}^{(2)}$\tabularnewline
\hline 
$\left\lceil -2\log(H)\right\rceil =5$ & $\cfrac{1}{10}$ & $\cfrac{1}{10}$ & 3.83\% & 6.24\%\tabularnewline
\hline 
$\left\lceil -2\log(H)\right\rceil =6$ & $\cfrac{1}{20}$ & $\cfrac{1}{20}$ & 2.20\% & 2.50\%\tabularnewline
\hline 
$\left\lceil -2\log(H)\right\rceil =8$ & $\cfrac{1}{40}$ & $\cfrac{1}{40}$ & 0.88\% & 0.90\%\tabularnewline
\hline 
\end{tabular}
\caption{Error comparison for Case 1 with a fixed contrast.}
\label{tab:case1_fixed}
\end{table}

In Table \ref{tab:case1_eff}, we present the results for effective properties
that are computed. First, we note that the scalings of these quantities
are in accordance with our theoretical findings. In Case 1, we
observe anisotropy as expected. $\alpha_{11}^{mn}\approx 0$, unless
$i=j=2$ since the flow in vertical direction for the gradient.
We observe larger $\alpha_{22}^{22}$ since the flow in the vertical direction
and the second continuum accounts for high conductivity. 
 We only show $\beta_{11}$ as other
$\beta_{ij}$'s depend on them and are similar
 (follows from symmetry and the fact that
the sum of elements in each row is zero). 
We note that with our scalings of conductivity, $\beta$ should scale
as $1/\epsilon$, which we observe in the table.

\begin{table}
\centering
\begin{tabular}{|c|c|c|c|c|c|}
\hline 
$l$ & $H$ & $\epsilon$ & $\alpha_{11}^{11}/|R_\omega|$ & $\alpha_{11}^{21}/|R_\omega|$ & $\alpha_{11}^{22}/|R_\omega|$\tabularnewline
\hline 
$\left\lceil -2\log(H)\right\rceil =5$ & $\cfrac{1}{10}$ & $\cfrac{1}{10}$ & $\approx0$ & $\approx0$ & 1.9685e-05\tabularnewline
\hline 
$\left\lceil -2\log(H)\right\rceil =5$ & $\cfrac{1}{10}$ & $\cfrac{1}{20}$ & $\approx0$ & $\approx0$ & 1.9006e-05\tabularnewline
\hline 
$\left\lceil -2\log(H)\right\rceil =5$ & $\cfrac{1}{10}$ & $\cfrac{1}{40}$ & $\approx0$ & $\approx0$ & 1.9247e-05\tabularnewline
\hline 
\end{tabular}

\begin{tabular}{|c|c|c|c|c|c|}
\hline 
$l$ & $H$ & $\epsilon$ & $\alpha_{22}^{11}/|R_\omega|$ & $\alpha_{22}^{21}/|R_\omega|$ & $\alpha_{22}^{22}/|R_\omega|$\tabularnewline
\hline 
$\left\lceil -2\log(H)\right\rceil =5$ & $\cfrac{1}{10}$ & $\cfrac{1}{10}$ & 1.4876e-05 & $\approx0$ & 0.0201\tabularnewline
\hline 
$\left\lceil -2\log(H)\right\rceil =5$ & $\cfrac{1}{10}$ & $\cfrac{1}{20}$ & 6.2843e-06 & $\approx0$ & 0.0401\tabularnewline
\hline 
$\left\lceil -2\log(H)\right\rceil =5$ & $\cfrac{1}{10}$ & $\cfrac{1}{40}$ & 3.1331e-06 & $\approx0$ & 0.0802\tabularnewline
\hline 
\end{tabular}

\begin{tabular}{|c|c|c|c|c|c|}
\hline 
$l$ & $H$ & $\epsilon$ & $\alpha_{12}^{11}/|R_\omega|$ & $\alpha_{12}^{12}/|R_\omega|$ & $\alpha_{12}^{22}/|R_\omega|$\tabularnewline
\hline 
$\left\lceil -2\log(H)\right\rceil =5$ & $\cfrac{1}{10}$ & $\cfrac{1}{10}$ & $\approx0$ & $\approx0$ & -1.1669e-05\tabularnewline
\hline 
$\left\lceil -2\log(H)\right\rceil =5$ & $\cfrac{1}{10}$ & $\cfrac{1}{20}$ & $\approx0$ & $\approx0$ & -1.4998e-05\tabularnewline
\hline 
$\left\lceil -2\log(H)\right\rceil =5$ & $\cfrac{1}{10}$ & $\cfrac{1}{40}$ & $\approx0$ & $\approx0$ & -1.7243e-05\tabularnewline
\hline 
\end{tabular}

\begin{tabular}{|c|c|c|c|c|c|}
\hline 
$l$ & $H$ & $\epsilon$ & $\beta_{11}/|R_\omega|$ & $\beta_{12}/|R_\omega|$ & $\beta_{22}/|R_\omega|$\tabularnewline
\hline 
$\left\lceil -2\log(H)\right\rceil =5$ & $\cfrac{1}{10}$ & $\cfrac{1}{10}$ & 0.0150 & -0.0150 & 0.0150\tabularnewline
\hline 
$\left\lceil -2\log(H)\right\rceil =5$ & $\cfrac{1}{10}$ & $\cfrac{1}{20}$ & 0.0301 & -0.0301 & 0.0301\tabularnewline
\hline 
$\left\lceil -2\log(H)\right\rceil =5$ & $\cfrac{1}{10}$ & $\cfrac{1}{40}$ & 0.0610 & -0.0610 & 0.0610\tabularnewline
\hline 
\end{tabular}
\caption{Effective properties $\alpha$ and $\beta$'s.  Case 1.}
\label{tab:case1_eff}
\end{table}

\subsection{Case 2}

\begin{figure}
\centering  
\includegraphics[scale=0.3]{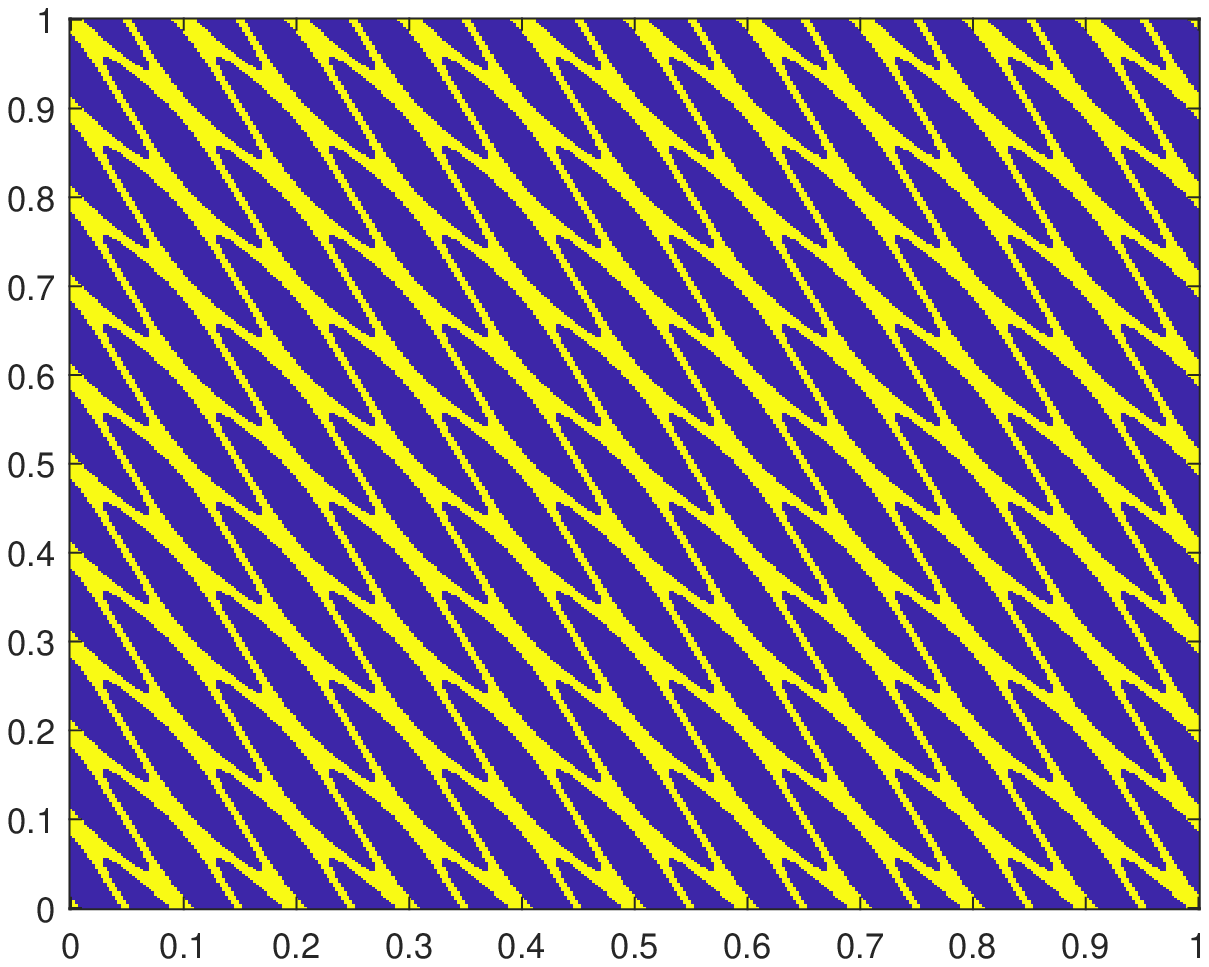}
\includegraphics[scale=0.3]{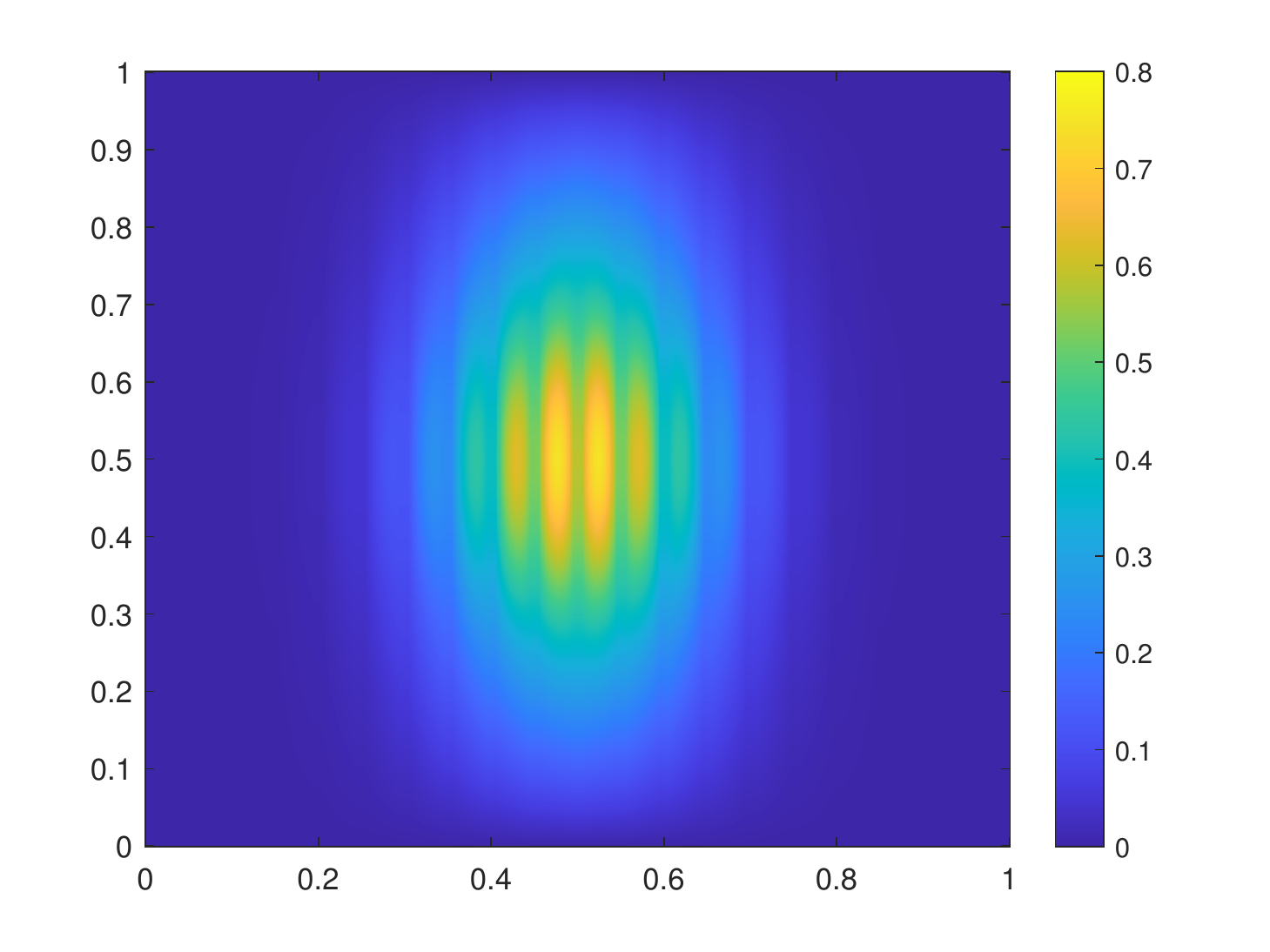}
\caption{Left: $\kappa$ for Case 2  with $\epsilon=\cfrac{1}{10}$. Right: The solution snapshot.}
\label{fig:case2}
\end{figure}

\begin{figure}
\centering
\includegraphics[scale=0.3]{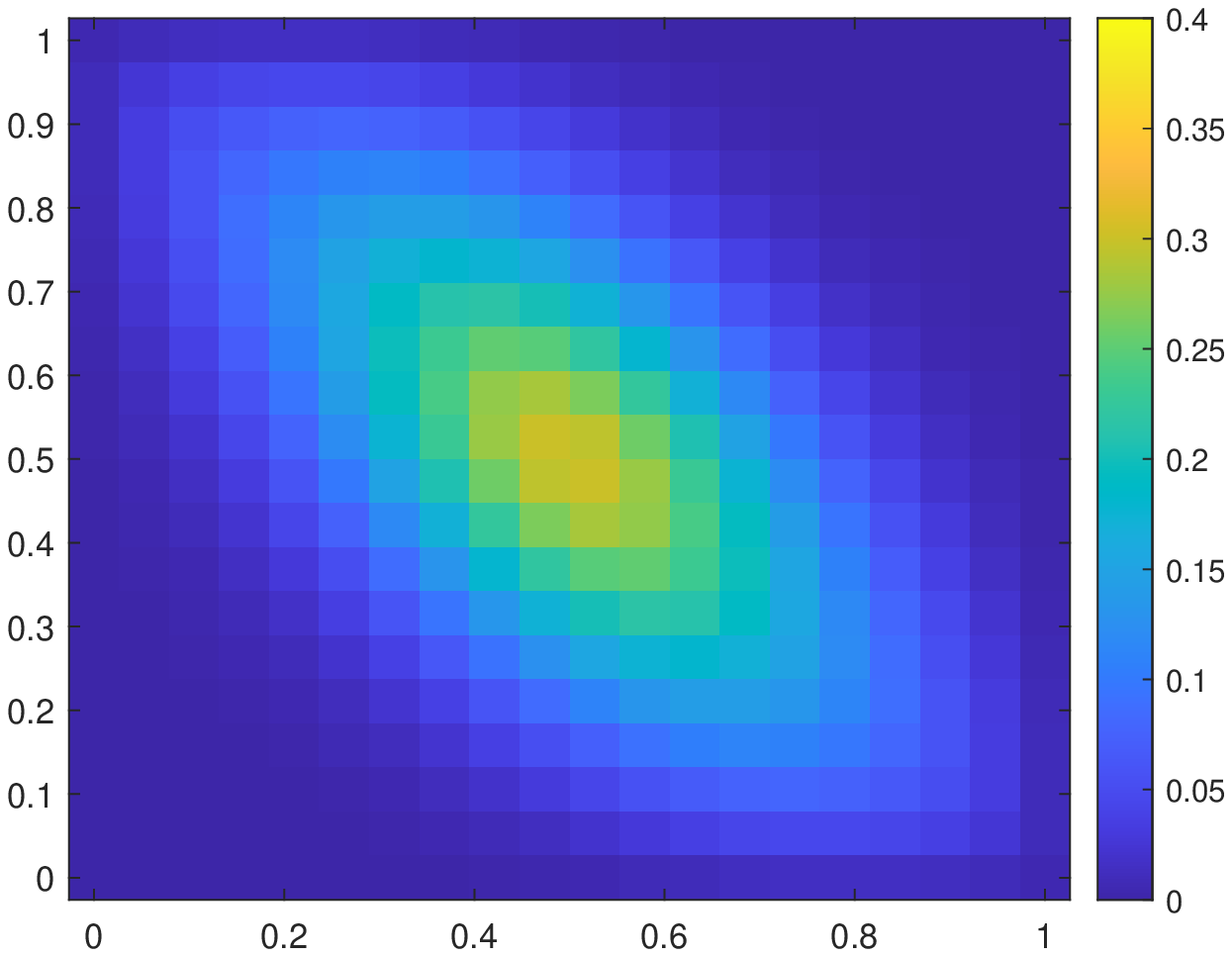} 
\includegraphics[scale=0.3]{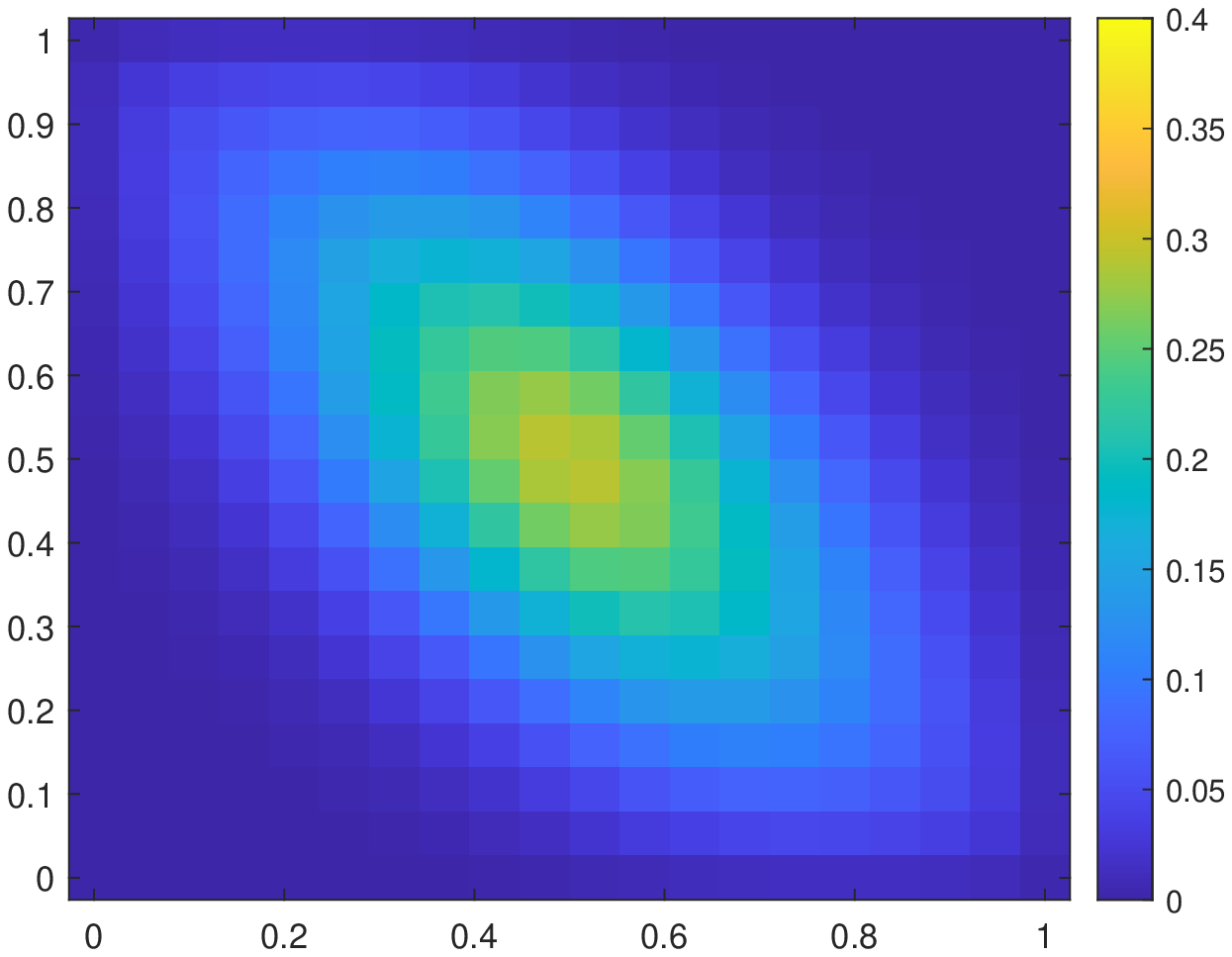}
\includegraphics[scale=0.3]{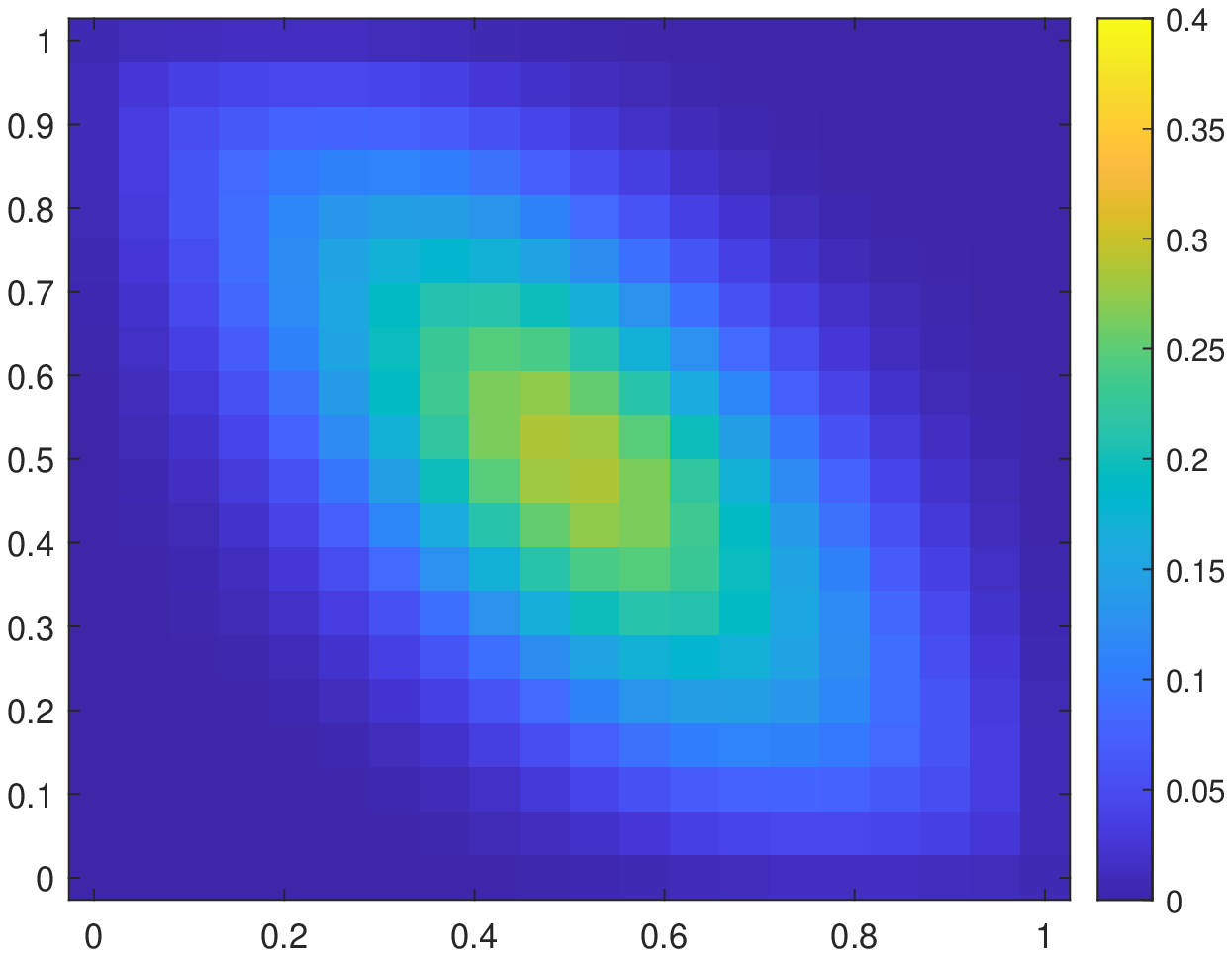} 
\includegraphics[scale=0.3]{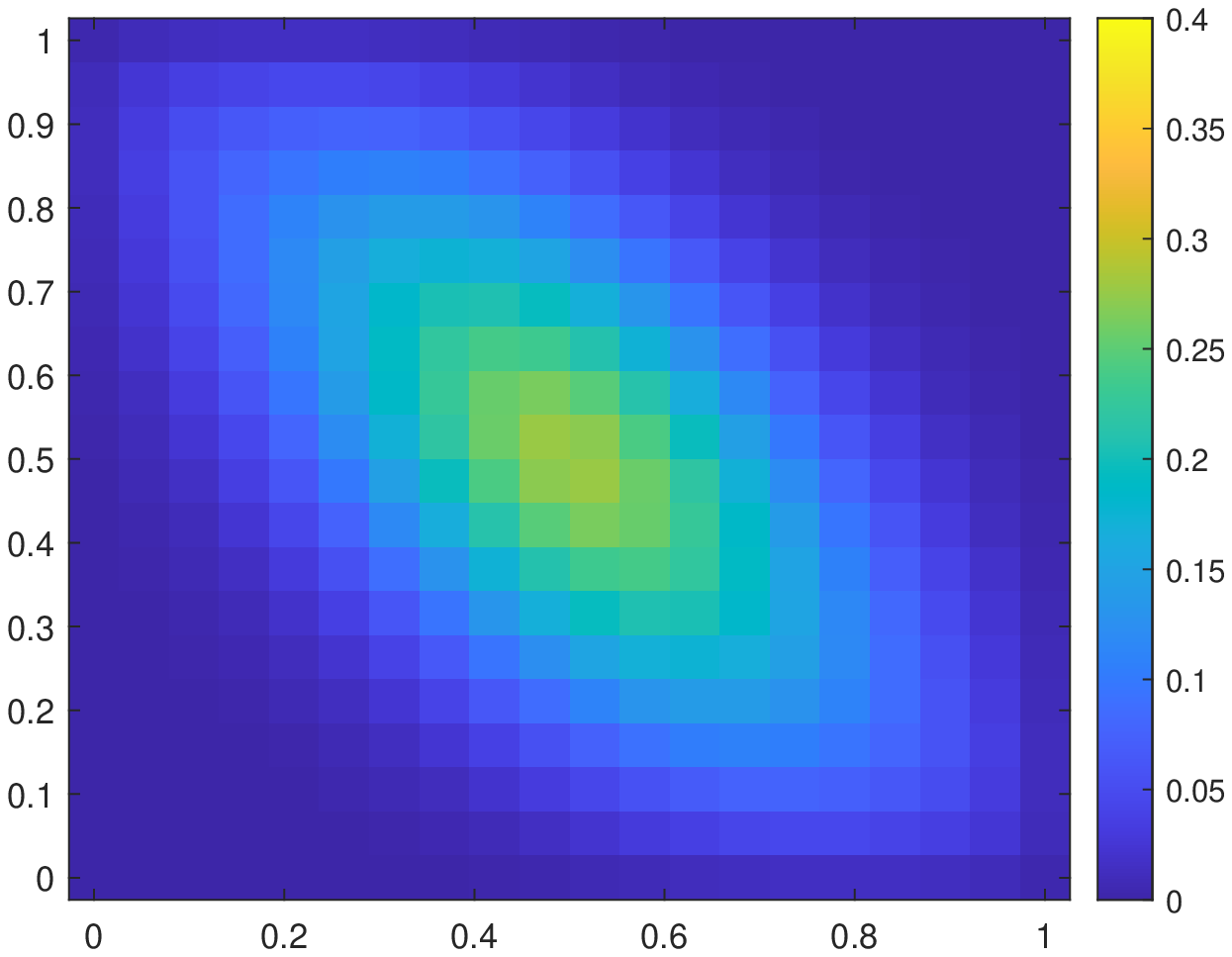}
\caption{ Case 2. Top-Left:
Reference averaged solution in $\Omega_{1}$. Top-Right: Multiscale
solution $U_{1}$. Bottom-Left: Reference averaged solution in $\Omega_{2}$.
Bottom-Right: Multiscale solution $U_{2}$.}
\label{fig:compare_case2}
\end{figure}

Next, we consider a different case, Case 2. We depict 
the conductivity field in Figure \ref{fig:case2}
and the corresponding fine-grid solution.
In Figure \ref{fig:compare_case2}, we depict upscaled solutions
and corresponding averaged fine-scale solutions. We observe that these
solutions are very close.
For this case, we present $e_2^{(i)}$ in Table \ref{tab:case2}.
We observe similar findings. 
First, we observe that the proposed
approach provides an accurate approximation of the averaged solution.
In this table, we decrease the coarse mesh size while keeping 
$\epsilon$ fixed. We observe that the error decreases if
we scale the number of layers as $-\log(H)$. In the second part
of the table, we decrease $\epsilon$ and observe that error does not change
if the number of layers is fixed. 
Finally, we reduce both $H$ and $\epsilon$ and observe that the error decreases.
Thus, the proposed
method is robust with respect to the size of heterogeneities.

\begin{table}
\centering
\begin{tabular}{|c|c|c|c|c|}
\hline 
$l$ & $H$ & $\epsilon$ & $e_{2}^{(1)}$ & $e_{2}^{(2)}$\tabularnewline
\hline 
$\left\lceil -2\log(H)\right\rceil =5$ & $\cfrac{1}{10}$ & $\cfrac{1}{40}$ & 8.42\% & 8.54\%\tabularnewline
\hline 
$\left\lceil -2\log(H)\right\rceil =6$ & $\cfrac{1}{20}$ & $\cfrac{1}{40}$ & 2.42\% & 2.50\%\tabularnewline
\hline 
$\left\lceil -2\log(H)\right\rceil =8$ & $\cfrac{1}{40}$ & $\cfrac{1}{40}$ & 0.65\% & 0.72\%\tabularnewline
\hline 
\end{tabular} %
\begin{tabular}{|c|c|c|c|c|}
\hline 
$l$ & $H$ & $\epsilon$ & $e_{2}^{(1)}$ & $e_{2}^{(2)}$\tabularnewline
\hline 
$5$ & $\cfrac{1}{10}$ & $\cfrac{1}{10}$ & 8.28\% & 10.08\%\tabularnewline
\hline 
$5$ & $\cfrac{1}{10}$ & $\cfrac{1}{20}$ & 9.60\% & 10.17\%\tabularnewline
\hline 
5 & $\cfrac{1}{10}$ & $\cfrac{1}{40}$ & 8.42\% & 8.54\%\tabularnewline
\hline 
\end{tabular}

\begin{tabular}{|c|c|c|c|c|}
\hline 
$l$ & $H$ & $\epsilon$ & $e_{2}^{(1)}$ & $e_{2}^{(2)}$\tabularnewline
\hline 
$\left\lceil -2\log(H)\right\rceil =5$ & $\cfrac{1}{10}$ & $\cfrac{1}{10}$ & 8.28\% & 10.08\%\tabularnewline
\hline 
$\left\lceil -2\log(H)\right\rceil =6$ & $\cfrac{1}{20}$ & $\cfrac{1}{20}$ & 2.98\% & 3.43\%\tabularnewline
\hline 
$\left\lceil -2\log(H)\right\rceil =8$ & $\cfrac{1}{40}$ & $\cfrac{1}{40}$ & 0.65\% & 0.72\%\tabularnewline
\hline 
$\left\lceil -2\log(H)\right\rceil =9$ & $\cfrac{1}{80}$ & $\cfrac{1}{80}$ & 0.18\% & 0.19\%\tabularnewline
\hline 
\end{tabular}
\caption{Error comparison for Case 2.}
\label{tab:case2}.
\end{table}

In Table \ref{tab:case2_eff}, numerical results for effective properties
are presented.  
In this case, we do not have strong anistropy and observe similar 
values for $\alpha$'s when the continua are fixed. From the second table,
we observe that the values of $\alpha_{22}^{mn}$ are larger compared to
$\alpha_{11}^{mn}$ and $\alpha_{12}^{mn}$. This is because the second continuum
account for high conductivity region.
Again,
 $\beta$ should scale
as $1/\epsilon$, which we observe in the table.

\begin{table}
\centering
\begin{tabular}{|c|c|c|c|c|c|}
\hline 
$l$ & $H$ & $\epsilon$ & $\alpha_{11}^{11}/|R_\omega|$ & $\alpha_{11}^{21}/|R_\omega|$ & $\alpha_{11}^{22}/|R_\omega|$\tabularnewline
\hline 
$\left\lceil -2\log(H)\right\rceil =5$ & $\cfrac{1}{10}$ & $\cfrac{1}{10}$ & 1.4528e-04 & -1.0341e-06 & 7.5875e-05\tabularnewline
\hline 
$\left\lceil -2\log(H)\right\rceil =5$ & $\cfrac{1}{10}$ & $\cfrac{1}{20}$ & 1.4534e-04 & -4.9526e-07 & 7.5874e-05\tabularnewline
\hline 
$\left\lceil -2\log(H)\right\rceil =5$ & $\cfrac{1}{10}$ & $\cfrac{1}{40}$ & 1.4535e-04 & -2.4487e-07 & 7.5866e-05\tabularnewline
\hline 
\end{tabular}

\begin{tabular}{|c|c|c|c|c|c|}
\hline 
$l$ & $H$ & $\epsilon$ & $\alpha_{22}^{11}/|R_\omega|$ & $\alpha_{22}^{12}/|R_\omega|$ & $\alpha_{22}^{22}/|R_\omega|$\tabularnewline
\hline 
$\left\lceil -2\log(H)\right\rceil =5$ & $\cfrac{1}{10}$ & $\cfrac{1}{10}$ & 0.0094 & -0.0119 & 0.0192\tabularnewline
\hline 
$\left\lceil -2\log(H)\right\rceil =5$ & $\cfrac{1}{10}$ & $\cfrac{1}{20}$ & 0.0187 & -0.0239 & 0.0384\tabularnewline
\hline 
$\left\lceil -2\log(H)\right\rceil =5$ & $\cfrac{1}{10}$ & $\cfrac{1}{40}$ & 0.0373 & -0.0477 & 0.0768\tabularnewline
\hline 
\end{tabular}

\begin{tabular}{|c|c|c|c|c|c|}
\hline 
$l$ & $H$ & $\epsilon$ & $\alpha_{12}^{11}/|R_\omega|$ & $\alpha_{12}^{12}/|R_\omega|$ & $\alpha_{12}^{22}/|R_\omega|$\tabularnewline
\hline 
$\left\lceil -2\log(H)\right\rceil =5$ & $\cfrac{1}{10}$ & $\cfrac{1}{10}$ & -1.4056e-04 & 2.5701e-05 & -5.7296e-05 \tabularnewline
\hline 
$\left\lceil -2\log(H)\right\rceil =5$ & $\cfrac{1}{10}$ & $\cfrac{1}{20}$ & -1.4302e-04 & 1.2903e-05 & -6.6648e-05\tabularnewline
\hline 
$\left\lceil -2\log(H)\right\rceil =5$ & $\cfrac{1}{10}$ & $\cfrac{1}{40}$ & -1.4419e-04 & 6.4525e-06 & -7.1260e-05\tabularnewline
\hline 
\end{tabular}

\begin{tabular}{|c|c|c|c|c|c|}
\hline 
$l$ & $H$ & $\epsilon$ & $\beta_{11}/|R_\omega|$ & $\beta_{12}/|R_\omega|$ & $\beta_{22}/|R_\omega|$\tabularnewline
\hline 
$\left\lceil -2\log(H)\right\rceil =5$ & $\cfrac{1}{10}$ & $\cfrac{1}{10}$ & 0.1572 & -0.1572 & 0.1572\tabularnewline
\hline 
$\left\lceil -2\log(H)\right\rceil =5$ & $\cfrac{1}{10}$ & $\cfrac{1}{20}$ & 0.3145 & -0.3145 & 0.3145\tabularnewline
\hline 
$\left\lceil -2\log(H)\right\rceil =5$ & $\cfrac{1}{10}$ & $\cfrac{1}{40}$ & 0.6291 & -0.6291 & 0.6291\tabularnewline
\hline 
\end{tabular}
\caption{Effective properties $\alpha$ and $\beta$'s. Case 2.}
\label{tab:case2_eff}
\end{table}

\subsection{Case 3}

\begin{figure}
\centering
\includegraphics[scale=0.3]{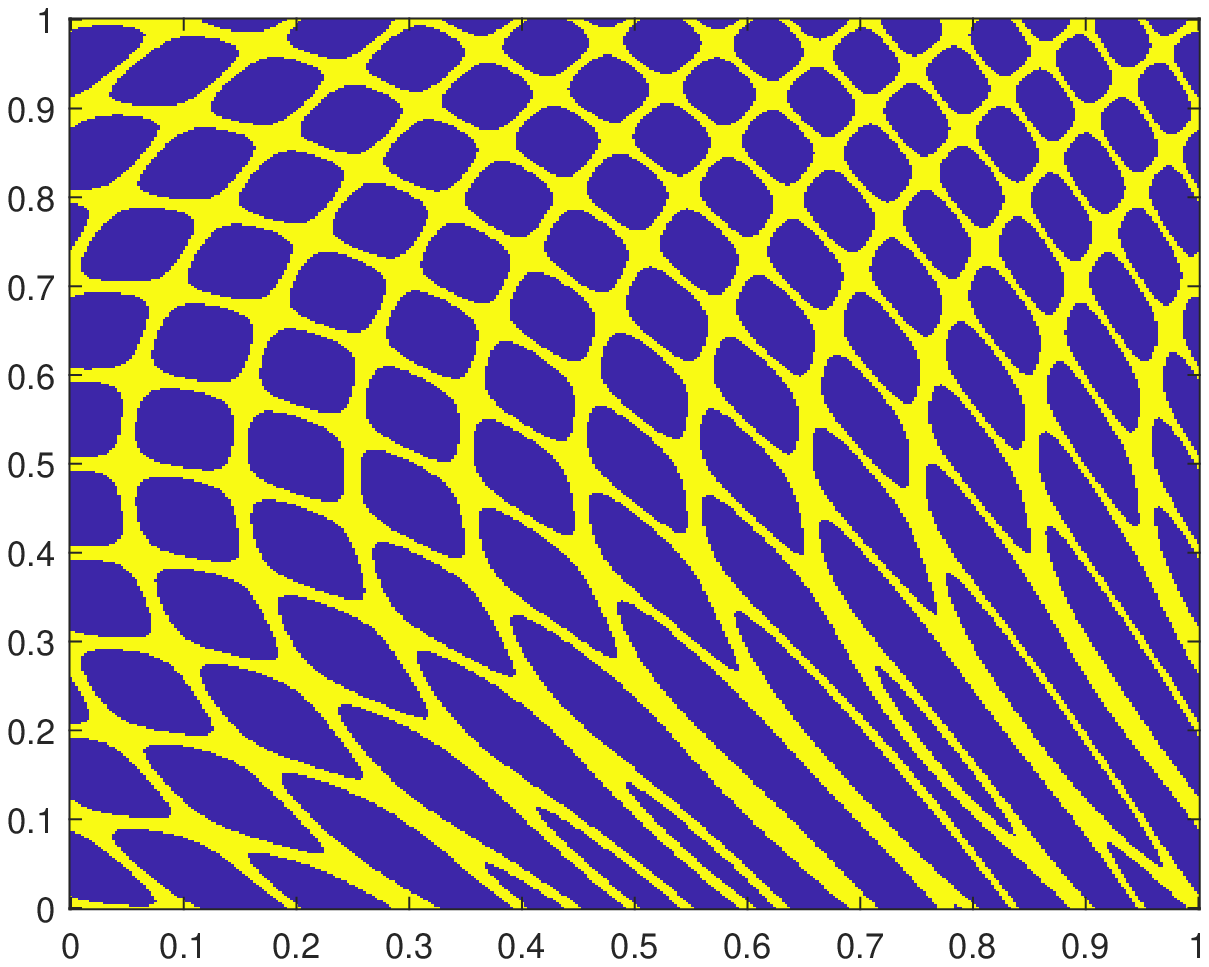}
\includegraphics[scale=0.3]{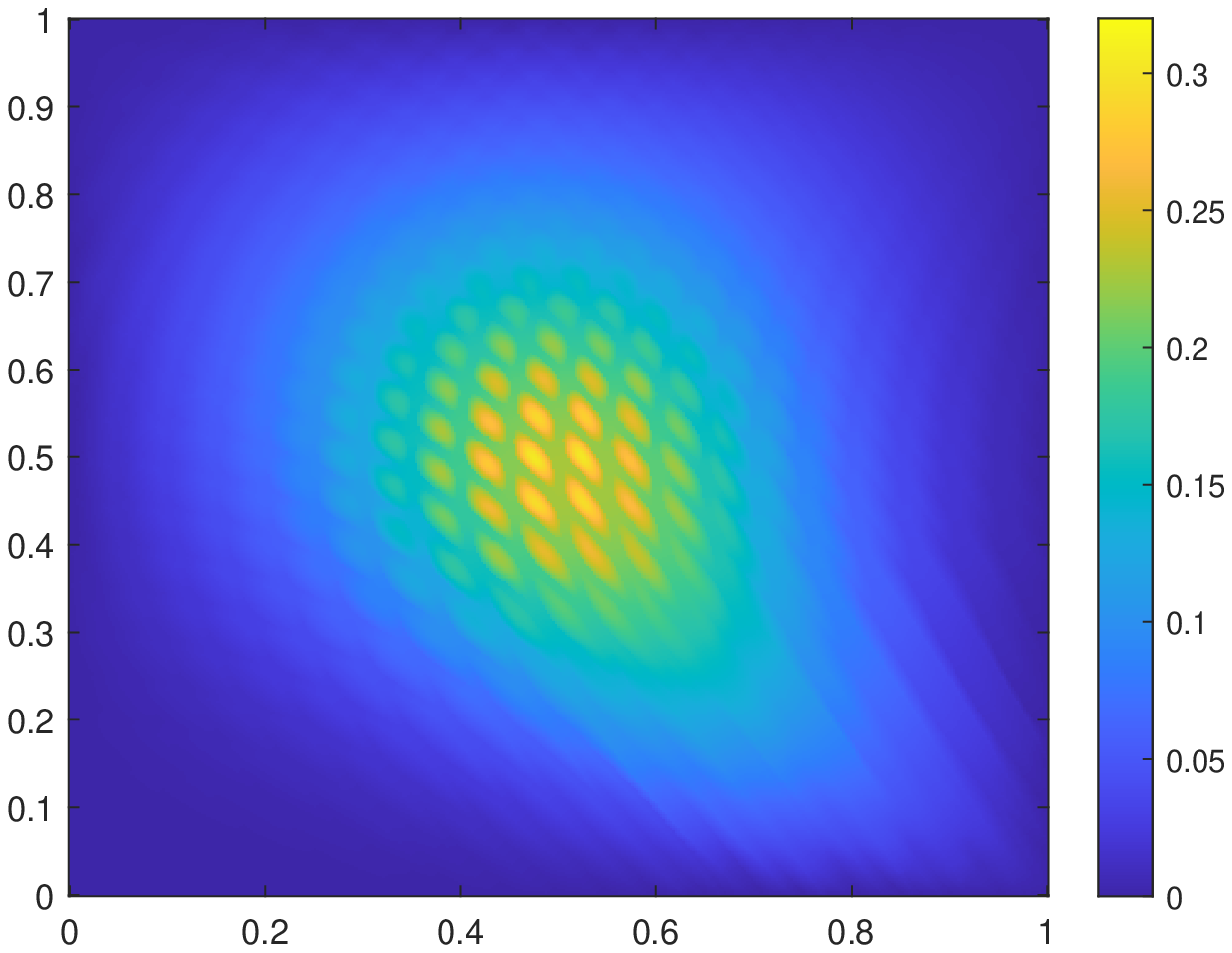}
\caption{Left: $\kappa$ for Case 3 with $\epsilon=\cfrac{1}{10}$. Right: The solution snapshot.}
\label{fig:case3}
\end{figure}

\begin{figure}
\centering
\includegraphics[scale=0.3]{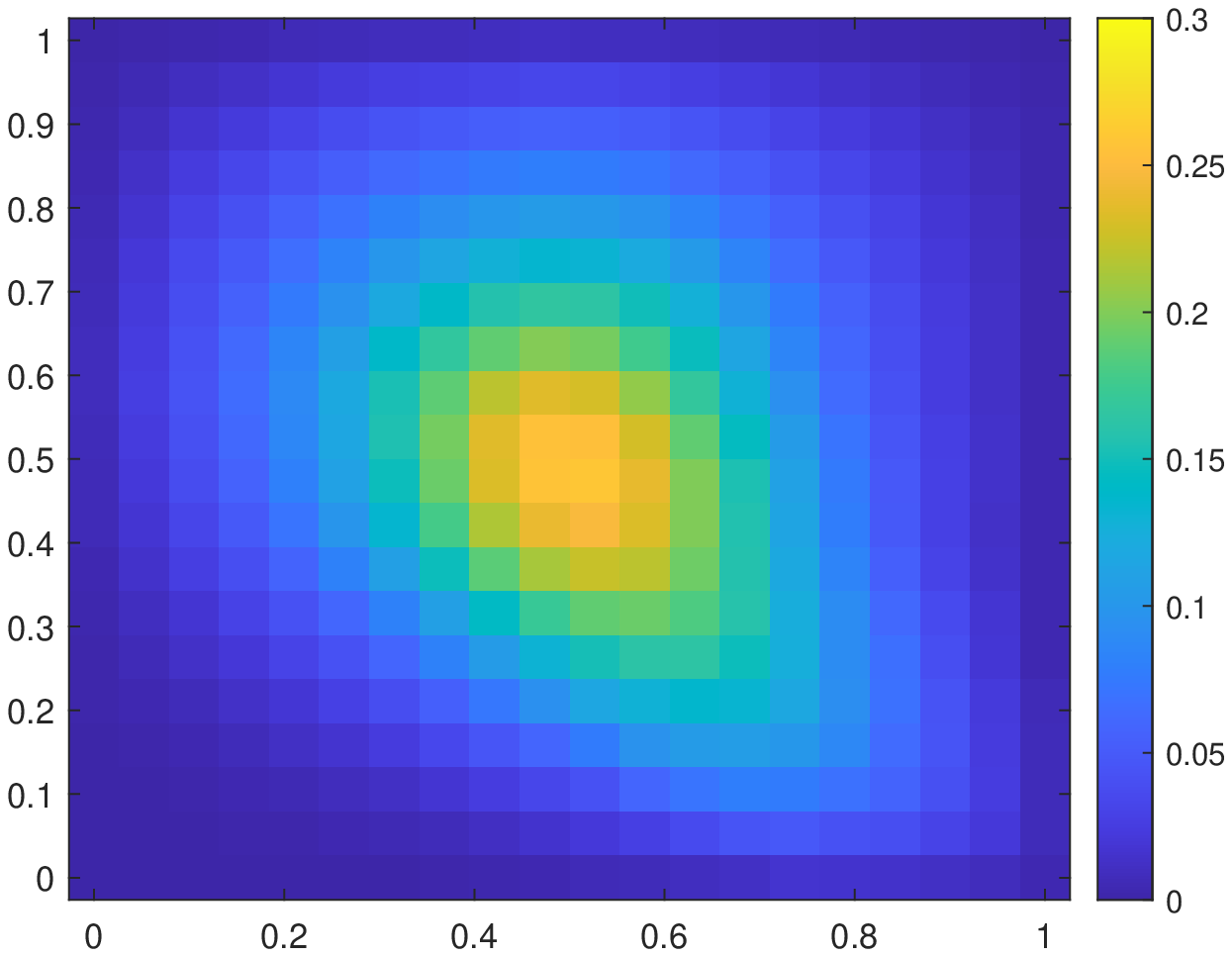} 
\includegraphics[scale=0.3]{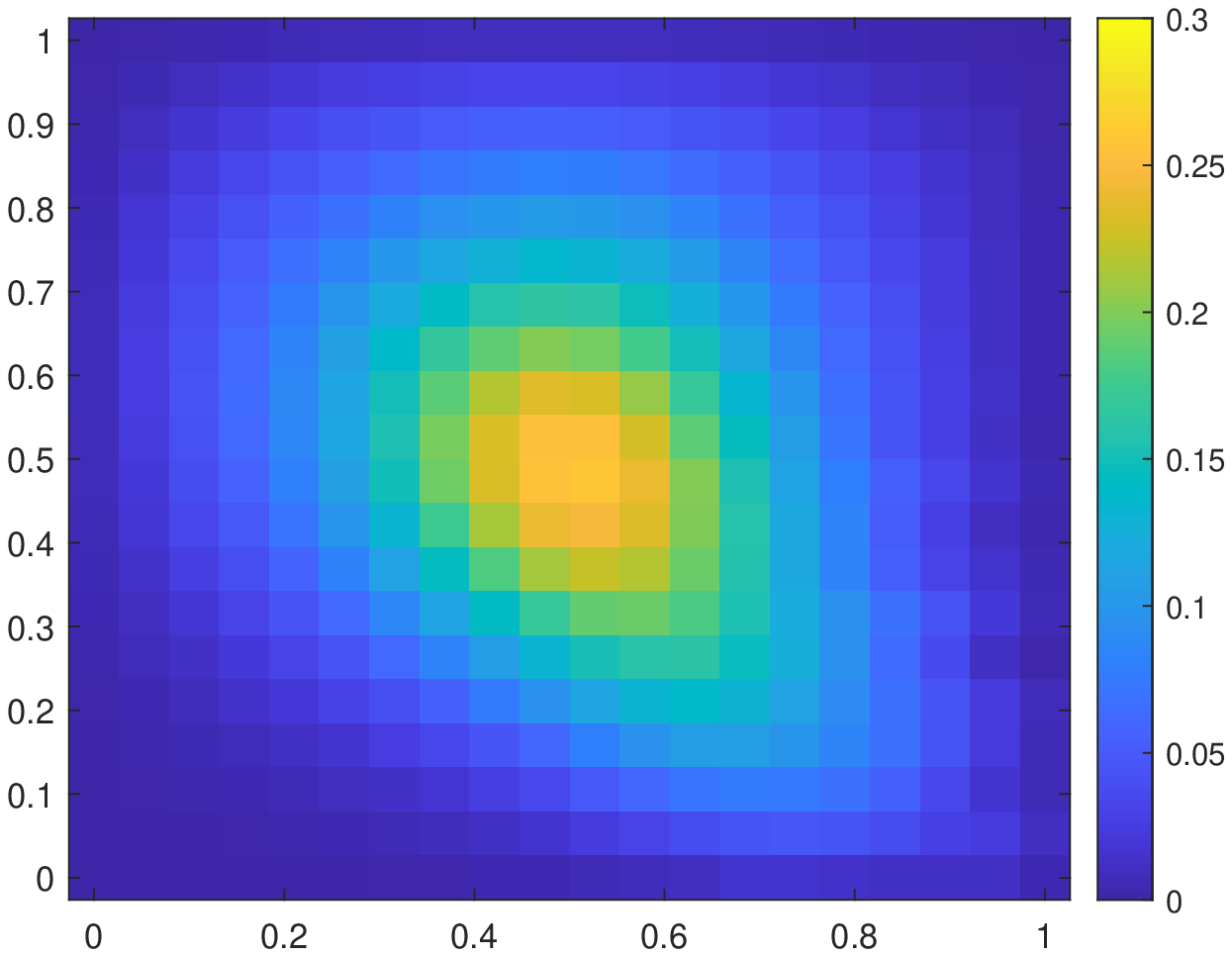}
\includegraphics[scale=0.3]{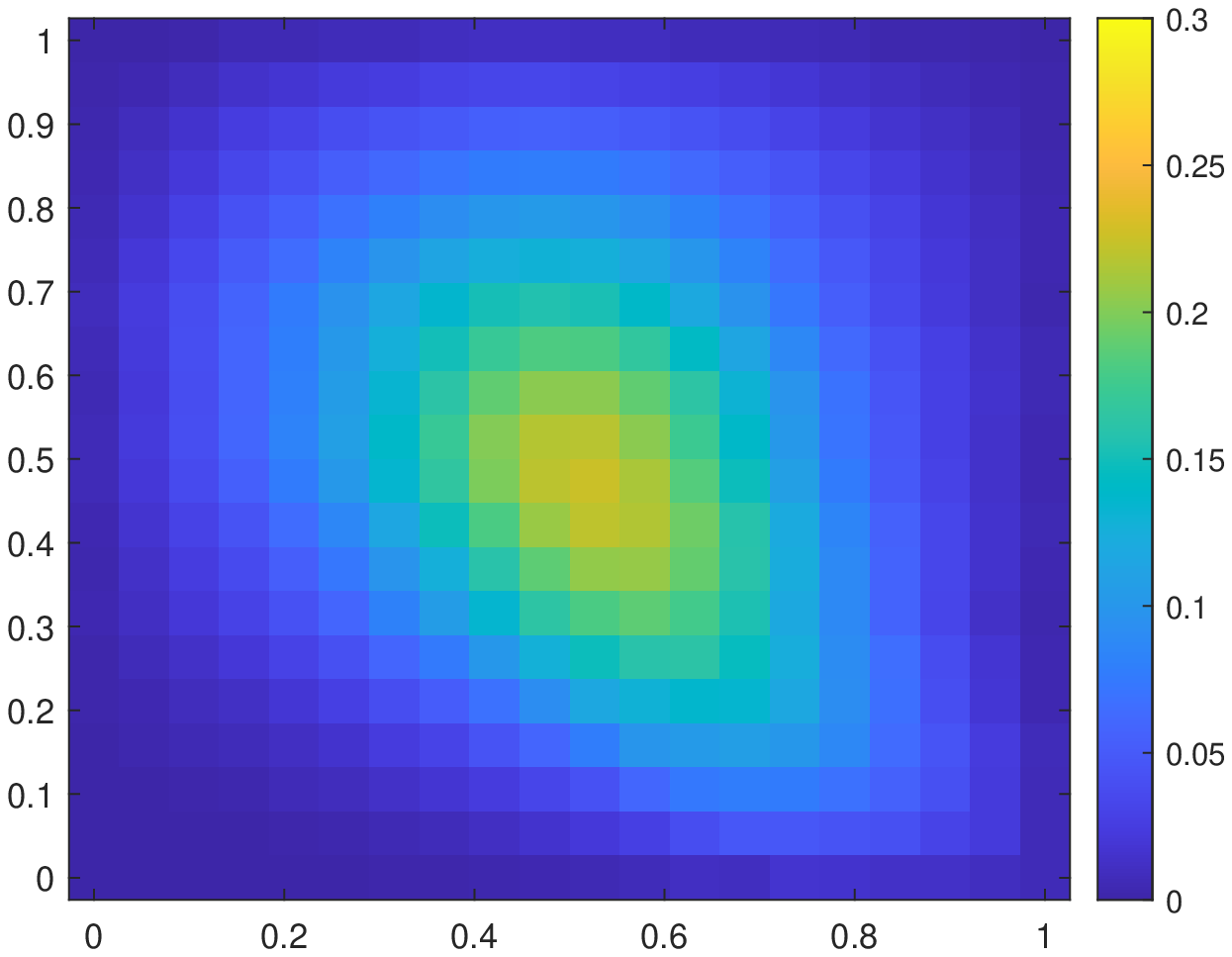}
 \includegraphics[scale=0.3]{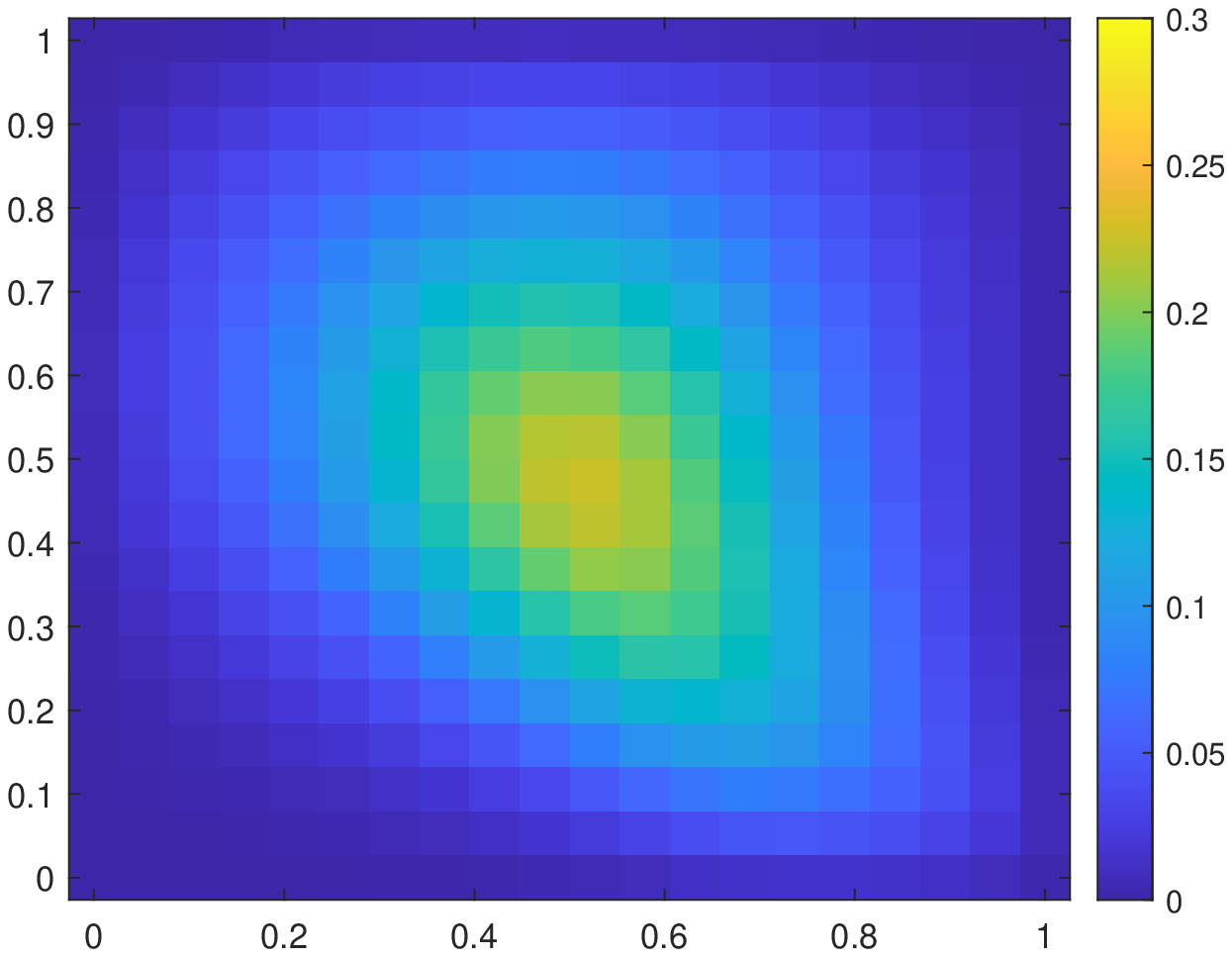}
\caption{Top-Left: Reference averaged solution in $\Omega_{1}$.
Top-Right: Multiscale
solution $U_{1}$. Bottom-Left: Reference averaged solution in $\Omega_{2}$.
Bottom-Right: Multiscale solution $U_{2}$.}
\label{fig:compare_case3}
\end{figure}

Next, we consider a different case, Case 3. We depict 
the conductivity field in Figure \ref{fig:case3}
and corresponding fine-grid solution. In this case, we do not 
have strict periodicity as in other cases and the conductivity
slowly changes.
In Figure \ref{fig:compare_case3}, we depict upscaled solutions
and corresponding averaged fine-scale solutions. We observe that these
solutions are approximately the same.
For this case, we present $e_2^{(i)}$ in Table \ref{tab:case3}.
We observe similar findings as in previous cases.
First, we observe that the proposed
approach provides an accurate approximation of the averaged solution.
In this table, we decrease the coarse mesh size while keeping 
$\epsilon$ fixed. We observe that the error decreases if
we scale the number of layers as $-\log(H)$. In the second part
of the table, we decrease $\epsilon$ and observe that error does not change
if the number of layers is fixed. 
Finally, we reduce both $H$ and $\epsilon$ and observe that the error decreases.
Thus, the proposed
method is robust with respect to the size of heterogeneities.

\begin{table}
\centering

\begin{tabular}{|c|c|c|c|c|}
\hline 
$l$ & $H$ & $\epsilon$ & $e_{2}^{(1)}$ & $e_{2}^{(2)}$\tabularnewline
\hline 
$\left\lceil -2\log(H)\right\rceil =5$ & $\cfrac{1}{10}$ & $\cfrac{1}{40}$ & 4.19\% & 4.24\%\tabularnewline
\hline 
$\left\lceil -2\log(H)\right\rceil =6$ & $\cfrac{1}{20}$ & $\cfrac{1}{40}$ & 1.21\% & 1.37\%\tabularnewline
\hline 
$\left\lceil -2\log(H)\right\rceil =8$ & $\cfrac{1}{40}$ & $\cfrac{1}{40}$ & 0.95\% & 1.14\%\tabularnewline
\hline 
\end{tabular} %
\begin{tabular}{|c|c|c|c|c|}
\hline 
$l$ & $H$ & $\epsilon$ & $e_{2}^{(1)}$ & $e_{2}^{(2)}$\tabularnewline
\hline 
$5$ & $\cfrac{1}{10}$ & $\cfrac{1}{10}$ & 5.41\% & 5.95\%\tabularnewline
\hline 
$5$ & $\cfrac{1}{10}$ & $\cfrac{1}{20}$ & 4.07\% & 4.21\%\tabularnewline
\hline 
5 & $\cfrac{1}{10}$ & $\cfrac{1}{40}$ & 4.19\% & 4.24\%\tabularnewline
\hline 
\end{tabular}

\begin{tabular}{|c|c|c|c|c|}
\hline 
$l$ & $H$ & $\epsilon$ & $e_{2}^{(1)}$ & $e_{2}^{(2)}$\tabularnewline
\hline 
$\left\lceil -2\log(H)\right\rceil =5$ & $\cfrac{1}{10}$ & $\cfrac{1}{10}$ & 5.41\% & 5.95\%\tabularnewline
\hline 
$\left\lceil -2\log(H)\right\rceil =6$ & $\cfrac{1}{20}$ & $\cfrac{1}{20}$ & 1.92\% & 2.25\%\tabularnewline
\hline 
$\left\lceil -2\log(H)\right\rceil =8$ & $\cfrac{1}{40}$ & $\cfrac{1}{40}$ & 0.95\% & 1.14\%\tabularnewline
\hline 
\end{tabular}
\caption{Error comparison for Case 3.}
\label{tab:case3}
\end{table}

\section{Conclusions}
In this paper, we propose a derivation of multicontinuum models using 
constraint cell problems in oversampled regions. The proposed cell problems 
allow reducing boundary effects and take into account both average
and gradient constraints. Imposing constraints on averages allows a fast
decay of artificial boundary effects. Our derivations show that one obtains
coupled equations and derives the formula for exchange between continua.
The exchange terms in the form of the reaction scale as the square of
the inverse of
RVE size and, thus, they dominate. As a result, the solutions in these continua are equal unless the diffusive terms can balance the reaction terms. 
This occurs if the media have high contrast. We discuss these issues
and how one can use spectral problems to define the continua. 
Based on obtained multicontinuum models, we show that one needs 
high contrast to have different average values. 
We derive multicontinuum models for dynamic problems with dynamic cell
problems and discuss
nonlinear cases. In addition, we discuss the use of both average and gradient
constraints at the same time and its disadvantages.
We also briefly discuss nonlinear multicontinuum models.
 Numerical results are presented. Our numerical results show that the proposed approach provides an accurate representation of the solution's averages and converges as we decrease the mesh size. We study various parameter regimes and their influences on effective properties. Finally, we would like to note that the proposed approaches can be extended to problems without scale separation. These problems are carefully studied in our previous works \cite{chung2018non,chung2018constraint}.

\bibliographystyle{abbrv}
\bibliography{references,references4,references1,references2,references3,decSol}

\end{document}